\documentclass[A4j,11pt]{article}
\usepackage{amssymb,amsmath,amsthm,amscd,amsfonts}
\usepackage{graphicx}
\usepackage{graphics}

\begin{document}
\title{\textbf{Gromov-Hausdorff-like distance function\\
on the Moduli space of Riemannian manifolds\\
and Hamilton convergence
}}
\author{Naoyuki Koike}
\date{\today}
\maketitle

%
%
%

\begin{abstract}
In this paper, we discuss how a Gromov-Hausdorff-like distance function over the space of all isometric classes 
of compact $C^k$-Riemannian manifolds should be defined in the aspect of the Riemannan submanifold theory, 
where $k\geq 1$.  The most important fact in this discussion is as follows.  
The Hausdorff distance function between two spheres of mutually distinct radii isometrically embedded into 
the hypebolic space of curvature $c$ converges to zero as $c\to-\infty$.  
The key in the construction of the Gromov-Hausdorff-like distance function given in this paper is to define 
the distance of two $C^{k+1}$-isometric embeddings of distinct compact $C^k$-Riemannian manifolds into a higher 
dimensional Riemannian manifold by using the Hausdorff distance function in the tangent bundle of order $k+1$ 
equipped with the Sasaki metric.  Furthermore, we show that the convergence of a sequence of compact Riemannian 
manifolds with respect to this distance function coincides with the convergence in the sense of 
R. S. Hamilton.  
\end{abstract}

\vspace{0.1truecm}





\vspace{0.3truecm}

\section{Introduction} 
First we shall recall the Gromov-Hausdorff distance function introduced by M. Gromov (\cite{G1}, \cite{G2}).  
Denote by $\widetilde{\mathcal M}$ the set of all metric spaces and $\mathcal M_c$ the space of all isometric 
classes of compact metric spaces and $[(X,d)]$ the isometric class of a compact metric space $(X,d)$.  
For metric spaces $(X,d)$ and $(\widetilde X,\widetilde d)$, denote by 
${\rm Emb}_{d.p.}((X,d),(\widetilde X,\widetilde d))$ the space of all distance-preserving embeddings of $(X,d)$ 
into $(\widetilde X,\widetilde d)$.  Let $\widetilde{\mathcal M}_c(\widetilde X,\widetilde d)$ be the set of 
all compact subsets of a metric space $(\widetilde X,\widetilde d)$.  
The {\it Hausdorff distance function} $d_{H,(\widetilde X,\widetilde d)}$ over 
$\widetilde{\mathcal M}_c(\widetilde X,\widetilde d)$ is defined by 
$$\begin{array}{r}
d_{H,(Y,d_Y)}(K_1,K_2):=\inf\{\varepsilon>0\,|\,K_2\subset B(K_1,\varepsilon)\,\,\&\,\,
K_1\subset B(K_2,\varepsilon)\}\\
(K_1,K_2\in\widetilde{\mathcal M}_c((\widetilde X,\widetilde d))),
\end{array}$$
where $B(K_i,\varepsilon)$ denotes the $\varepsilon$-neighborhood of $K_i$.  
By using this distance function, the {\it Gromov-Hausdorff distance funcion} $d_{GH}$ over $\mathcal M_c$ 
is defined by 
{\small 
\begin{align*}
&d_{GH}([(X_1,d_1)],[(X_2,d_2)])\\
:=&\mathop{\inf}_{(\widetilde X,\widetilde d)\in\widetilde{\mathcal M}}\,
\inf\{d_{H,(\widetilde X,\widetilde d)}(f_1(X_1),f_2(X_2))\,|\,f_i\in{\rm Emb}_{d.p.}((X_i,d_i),
(\widetilde X,\widetilde d))\,\,\,(i=1,2)\},\end{align*}
}

\noindent
where $\inf\{d_{H,(\widetilde X,\widetilde d)}(f_1(X_1),f_2(X_2))\,|\,f_i\in{\rm Emb}_{d.p.}((X_i,d_i),
(\widetilde X,\widetilde d))\,\,(i=1,2)\}$ implies $\infty$ 
in the case where ${\rm Emb}_{d.p.}((X_1,d_1),(\widetilde X,\widetilde d))=\emptyset$ or 
${\rm Emb}_{d.p.}((X_2,d_2),(\widetilde X,\widetilde d))=\emptyset$.  
It is well-known that this function $d_{GH}$ gives a distance function over $\mathcal M_c$ and furthermore 
$(\mathcal M_c,d_{GH})$ is a complete metric space.  

In this paper, we introduce a Gromov-Hausdorff-like distance function over the space of 
all isometric classes of compact $C^k$-Riemannian manifolds in the aspect of the Riemannian submanifold theory, 
where $k\geq 1$.  Furthermore, the convergence of compact Riemannian manifolds with respect to this distance 
function coincides with the convergence in the sense of R. S. Hamilton.  

\section{Some important examples} 
Let $k\in\mathbb N$ or $k=\infty$.  
Denote by $\widetilde{\mathcal RM}^k$ the set of all $C^k$-Riemannian manifolds and 
$\mathcal{RM}^k_c$ the space of all isometric classes of compact $C^k$-Riemannian manifolds and 
$[(M,g)]$ the isometric class of a compact Riemannian manifold $(M,g)$.  
Here we note that ``$C^k$-Riemannian manifold'' means a $C^{k+1}$-manifold equipped with a $C^k$-Riemannian metric.  
For $C^k$-Riemannian manifolds $(M,g)$ and $(\widetilde M,\widetilde g)$, denote by 
${\rm Emb}^{k+1}_I((M,g),(\widetilde M,\widetilde g))$ (resp. 
${\rm Emb}_I^{t.g.}((M,g),(\widetilde M,\widetilde g))$) 
the space of all $C^{k+1}$-isometric embeddings (resp. all totally geodesic $C^{k+1}$-isometric embeddings) of 
$(M,g)$ into $(\widetilde M,\widetilde g)$.  
Here we note that the following facts hold:

\vspace{0.5truecm}

\noindent
{\bf Fact.} {\sl Let $f$ be a $C^{k+1}$-isometric embedding of $(M,g)$ into $(\widetilde M,\widetilde g)$.  
If $f$ is not totally geodesic, then it is not a distance-preserving embedding of $(M,d_g)$ into 
$(\widetilde M,d_{\widetilde g})$, where $d_g$ (resp. $d_{\widetilde g}$) denotes the Riemannian distance function 
of $g$ (resp. $\widetilde g$).  
Even if $f$ is a totally geodesic $C^{k+1}$-isometric embedding, it is not necessarily 
a distance-preserving embedding of $(M,d_g)$ into $(\widetilde M,d_{\widetilde g})$ (see Figure 1).  
On the othe hand, if $f$ is a totally geodesic $C^{k+1}$-isometric embedding and if 
$(\widetilde M,\widetilde g)$ is a Hadamard manifold, then it is a distance-perserving embedding of $(M,d_g)$ into 
$(\widetilde M,d_{\widetilde g})$.}

\vspace{0.25truecm}

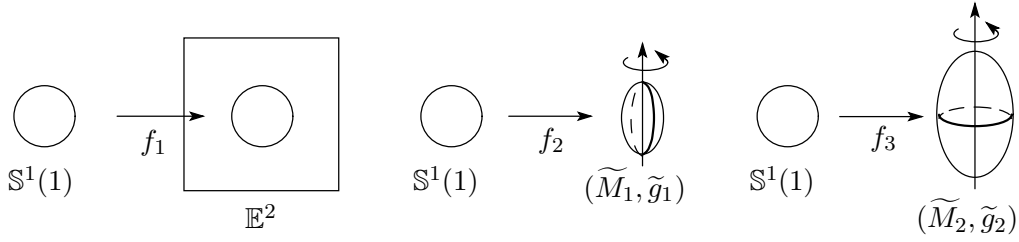
\centerline{
\unitlength 0.1in
\begin{picture}( 52.6000, 20.4000)(  8.8000,-28.6000)
%
\special{pn 8}%
\special{pa 1800 1000}%
\special{pa 2610 1000}%
\special{pa 2610 1800}%
\special{pa 1800 1800}%
\special{pa 1800 1000}%
\special{fp}%
%
\special{pn 8}%
\special{pa 1450 1410}%
\special{pa 1900 1410}%
\special{fp}%
\special{sh 1}%
\special{pa 1900 1410}%
\special{pa 1834 1390}%
\special{pa 1848 1410}%
\special{pa 1834 1430}%
\special{pa 1900 1410}%
\special{fp}%
\put(15.7000,-14.8000){\makebox(0,0)[lt]{$f_1$}}%
\put(21.2000,-18.7000){\makebox(0,0)[lt]{$\mathbb E^2$}}%
%
\special{pn 13}%
\special{pa 4200 1230}%
\special{pa 4226 1248}%
\special{pa 4240 1276}%
\special{pa 4250 1308}%
\special{pa 4254 1338}%
\special{pa 4258 1370}%
\special{pa 4260 1402}%
\special{pa 4260 1434}%
\special{pa 4258 1466}%
\special{pa 4254 1498}%
\special{pa 4250 1530}%
\special{pa 4240 1560}%
\special{pa 4228 1590}%
\special{pa 4204 1610}%
\special{pa 4200 1610}%
\special{sp}%
%
\special{pn 8}%
\special{pa 4202 1230}%
\special{pa 4234 1238}%
\special{pa 4258 1258}%
\special{pa 4278 1284}%
\special{pa 4292 1314}%
\special{pa 4302 1344}%
\special{pa 4308 1376}%
\special{pa 4310 1406}%
\special{pa 4310 1438}%
\special{pa 4306 1470}%
\special{pa 4300 1502}%
\special{pa 4290 1532}%
\special{pa 4274 1560}%
\special{pa 4254 1586}%
\special{pa 4228 1604}%
\special{pa 4198 1610}%
\special{pa 4168 1602}%
\special{pa 4142 1582}%
\special{pa 4124 1556}%
\special{pa 4110 1526}%
\special{pa 4100 1496}%
\special{pa 4094 1464}%
\special{pa 4090 1434}%
\special{pa 4090 1402}%
\special{pa 4094 1370}%
\special{pa 4102 1338}%
\special{pa 4112 1308}%
\special{pa 4126 1280}%
\special{pa 4146 1254}%
\special{pa 4172 1236}%
\special{pa 4200 1230}%
\special{sp}%
%
\special{pn 8}%
\special{pa 4200 1620}%
\special{pa 4176 1602}%
\special{pa 4162 1574}%
\special{pa 4152 1544}%
\special{pa 4146 1512}%
\special{pa 4142 1480}%
\special{pa 4140 1448}%
\special{pa 4140 1416}%
\special{pa 4142 1384}%
\special{pa 4146 1352}%
\special{pa 4152 1322}%
\special{pa 4160 1290}%
\special{pa 4174 1262}%
\special{pa 4196 1242}%
\special{pa 4200 1240}%
\special{sp 0.070}%
%
\special{pn 8}%
\special{pa 3500 1410}%
\special{pa 3930 1410}%
\special{fp}%
\special{sh 1}%
\special{pa 3930 1410}%
\special{pa 3864 1390}%
\special{pa 3878 1410}%
\special{pa 3864 1430}%
\special{pa 3930 1410}%
\special{fp}%
\put(36.5000,-14.7000){\makebox(0,0)[lt]{$f_2$}}%
%
\special{pn 13}%
\special{ar 5942 1402 190 60  6.2831853 6.2831853}%
\special{ar 5942 1402 190 60  0.0000000 3.1415927}%
%
\special{pn 8}%
\special{ar 5942 1422 190 60  3.1415927 3.6215927}%
\special{ar 5942 1422 190 60  3.9095927 4.3895927}%
\special{ar 5942 1422 190 60  4.6775927 5.1575927}%
\special{ar 5942 1422 190 60  5.4455927 5.9255927}%
\special{ar 5942 1422 190 60  6.2135927 6.2831853}%
%
\special{pn 8}%
\special{pa 5232 1412}%
\special{pa 5662 1412}%
\special{fp}%
\special{sh 1}%
\special{pa 5662 1412}%
\special{pa 5594 1392}%
\special{pa 5608 1412}%
\special{pa 5594 1432}%
\special{pa 5662 1412}%
\special{fp}%
\put(53.9000,-14.6000){\makebox(0,0)[lt]{$f_3$}}%
%
\special{pn 8}%
\special{ar 5940 1400 200 330  0.0000000 6.2831853}%
\put(8.8000,-16.7000){\makebox(0,0)[lt]{$\mathbb S^1(1)$}}%
\put(30.0000,-16.7000){\makebox(0,0)[lt]{$\mathbb S^1(1)$}}%
\put(47.6000,-16.7000){\makebox(0,0)[lt]{$\mathbb S^1(1)$}}%
\put(18.2000,-21.7000){\makebox(0,0)[lt]{{\small $f_1\in{\rm Emb}_I^{k+1}(\mathbb S^1(1),\mathbb E^2)\setminus{\rm Emb}_I^{t.g.}(\mathbb S^1(1),\mathbb E^2)$}}}%
\put(18.2000,-24.1000){\makebox(0,0)[lt]{{\small $f_2\in{\rm Emb}_I^{t.g.}(\mathbb S^1(1),(\widetilde M_1,\widetilde g_1))\setminus{\rm Emb}_{d.p.}(\mathbb S^1(1),(\widetilde M_1,\widetilde g_1))$}}}%
\put(18.2000,-26.4000){\makebox(0,0)[lt]{{\small $f_3\in{\rm Emb}_{d.p.}(\mathbb S^1(1),(\widetilde M_2,\widetilde g_2))$}}}%
\put(38.9000,-16.8000){\makebox(0,0)[lt]{$(\widetilde M_1,\widetilde g_1)$}}%
\put(56.3000,-18.4000){\makebox(0,0)[lt]{$(\widetilde M_2,\widetilde g_2)$}}%
\put(25.5000,-28.6000){\makebox(0,0)[lt]{$\,$}}%
%
\special{pn 8}%
\special{ar 3200 1410 160 160  0.0000000 6.2831853}%
%
\special{pn 8}%
\special{pa 4200 1670}%
\special{pa 4200 1030}%
\special{fp}%
\special{sh 1}%
\special{pa 4200 1030}%
\special{pa 4180 1098}%
\special{pa 4200 1084}%
\special{pa 4220 1098}%
\special{pa 4200 1030}%
\special{fp}%
%
\special{pn 8}%
\special{ar 4216 1128 116 42  5.7893611 6.2831853}%
\special{ar 4216 1128 116 42  0.0000000 3.7771380}%
%
\special{pn 8}%
\special{pa 4330 1120}%
\special{pa 4270 1070}%
\special{fp}%
\special{sh 1}%
\special{pa 4270 1070}%
\special{pa 4308 1128}%
\special{pa 4312 1104}%
\special{pa 4334 1098}%
\special{pa 4270 1070}%
\special{fp}%
%
\special{pn 8}%
\special{ar 5946 978 116 42  5.7893611 6.2831853}%
\special{ar 5946 978 116 42  0.0000000 3.7771380}%
%
\special{pn 8}%
\special{pa 6060 970}%
\special{pa 6000 920}%
\special{fp}%
\special{sh 1}%
\special{pa 6000 920}%
\special{pa 6038 978}%
\special{pa 6042 954}%
\special{pa 6064 948}%
\special{pa 6000 920}%
\special{fp}%
%
\special{pn 8}%
\special{pa 5940 1790}%
\special{pa 5940 820}%
\special{fp}%
\special{sh 1}%
\special{pa 5940 820}%
\special{pa 5920 888}%
\special{pa 5940 874}%
\special{pa 5960 888}%
\special{pa 5940 820}%
\special{fp}%
%
\special{pn 8}%
\special{ar 4960 1410 160 160  0.0000000 6.2831853}%
%
\special{pn 8}%
\special{ar 1070 1410 160 160  0.0000000 6.2831853}%
%
\special{pn 8}%
\special{ar 2210 1410 160 160  0.0000000 6.2831853}%
\end{picture}%
\hspace{0.25truecm}}

\vspace{0.2truecm}

\hspace{1truecm}{\bf Figure 1$\,\,:\,\,$ Gap between distance-preserving embeddings}

\hspace{3.3truecm}{\bf and isometric embeddings}

\vspace{0.25truecm}

\noindent
By using the Hausdorff distance functions, we define a function 
$\breve d^k_{GH}$ over $\mathcal{RM}^k_c\times\mathcal{RM}^k_c$ by 
{\small 
$$\begin{array}{l}
\hspace{0.5truecm}\displaystyle{\breve d^k_{GH}([(M_1,g_1)],[(M_2,g_2)])}\\
\displaystyle{:=\mathop{\inf}_{(\widetilde M,\widetilde g)\in\widetilde{\mathcal{RM}}^k}\,
\inf\{d_{H,(\widetilde M,d_{\widetilde g})}(f_1(M_1),f_2(M_2))\,|\,f_i\in{\rm Emb}^{k+1}_I((M_i,g_i),
(\widetilde M,\widetilde g))\,(i=1,2)\}.}
\end{array}$$}
This definition seems to be natural.  
However, we can show that $\breve d^k_{GH}$ is not a distance function over $\mathcal{RM}^k_c$.  
In fact, we can give some counter-examples as follows.  
Let $\mathbb S^n(r^{-2})$ be the sphere of radius $r$ centered at the origin $o=(0,\cdots,0)$ in the 
$(n+1)$-dimensional Euclidean space $\mathbb E^{n+1}$ (i.e., 
$\mathbb S^n(r^{-2})=\{(x_1,\cdots,x_{n+1})\in\mathbb E^{n+1}\,|\,x_1^2+\cdots+x_{n+1}^2=r^2\}$) and denote by 
$g_{\mathbb E}$ the Euclidean metric of $\mathbb E^{n+1}$.  
Let $\iota^S_r$ be the inclusion map of $\mathbb S^n(r^{-2})$ into $\mathbb E^{n+1}$ and denote by $g^S_r$ 
the induced metric $(\iota^S_r)^{\ast}g_{\mathbb E}$.  In the sequel, we abbreviate 
$(\mathbb S^n(r^{-2}),g^S_r)$ as $\mathbb S^n(r^{-2})$.  
In the case of $n\geq 2$, $\mathbb S^n(r^{-2})$ is the $n$-dimensional sphere of constant 
curvature $r^{-2}$.  We consider the case of $n=1$.  The length of $(\mathbb S^1(r^{-2}),g^S_r)$ is equal to 
$2\pi r$.  Fix positive numbers $r_1<r_2$.   It is clear that $\mathbb S^1(r_1^{-2})$ is not isometric to 
$\mathbb S^1(r_2^{-2})$.  However, we can show 
$$\inf\{d_{H,\mathbb E^2}(f_1(\mathbb S^1(r_1^{-2})),f_2(\mathbb S^1(r_2^{-2})))\,|\,
f_i\in{\rm Emb}^{k+1}_I(\mathbb S^1(r_i^{-2}),\mathbb E^2)\,(i=1,2)\}=0$$
as follows.  Take any positive number $\varepsilon$.  
Let $N_{\varepsilon}(\mathbb S^1(r_1^{-2}))$ be the $\varepsilon$-tubular neighborhood of 
$\mathbb S^1(r_1^{-2})$, that is, 
$$N_{\varepsilon}(\mathbb S^1(r_1^{-2}))
:=\{(x_1,x_2)\in\mathbb E^2\,|\,(r_1-\varepsilon)^2<x_1^2+x_2^2<(r_1+\varepsilon)^2\}.$$
For any sufficiently small postive number $\varepsilon$, $\mathbb S^1(r_2^{-2})$ can be isometrically embedded 
into $N_{\varepsilon}(\mathbb S^1(r_1^{-2}))$ (see Figure 2).  
This fact implies that 
$$\inf\{d_{H,\mathbb E^2}(f_1(\mathbb S^1(r_1^{-2})),f_2(\mathbb S^1(r_2^{-2})))\,|\,
f_i\in{\rm Emb}^{k+1}_I(\mathbb S^1(r_i^{-2}),\mathbb E^2)\,\,(i=1,2)\}=0.$$
Hence we obtain 
$$\breve d^k_{GH}([\mathbb S^1(r_1^{-2})],[\mathbb S^1(r_2^{-2})])=0.$$
Thus $\breve d^k_{GH}$ is not a distance function over $\mathcal{RM}^k_c$.  
For any Riemannian manifold $(M,g)$, we consider the product Riemannian manifolds 
$M\times \mathbb S^1(r_i^{-2})$ ($i=1,2$).  
In more general, it is shown that $[M\times \mathbb S^1(r_1^{-2})]\not=[M\times \mathbb S^1(r_2^{-2})]$ but 
$\breve d^k_{GH}([M\times \mathbb S^1(r_1^{-2})],[M\times \mathbb S^1(r_2^{-2})])=0$.  
Isometric embeddings $f_1$ and $f_2$ in Figure 2 are sufficiently close as $C^0$ embeddings but they are not 
close as $C^1$-embeddings and they are very far as $C^2$-embeddings.  On the other hand, isometric embeddings 
$\widehat f_1$ and $\widehat f_2$ in Figure 3 are sufficiently close as $C^{\infty}$-embeddings, where 
$2r_1<r_2<2r_1+\varepsilon$ ($\varepsilon:$ a sufficiently small positive number).  

\vspace{0.1truecm}

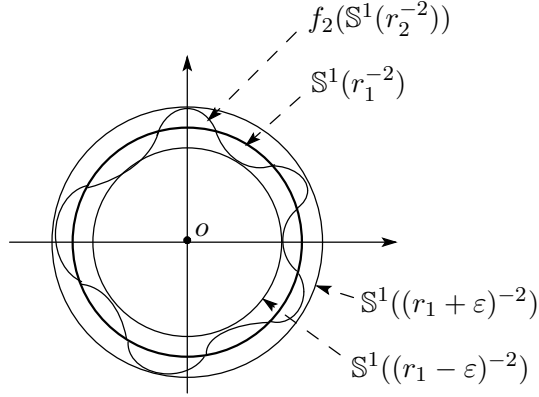
\centerline{
\unitlength 0.1in
\begin{picture}( 26.8000, 25.2000)(  8.1000,-28.7000)
\put(23.6000,-28.7000){\makebox(0,0)[lt]{$\,$}}%
%
\special{pn 13}%
\special{ar 2400 1600 600 600  0.0000000 6.2831853}%
%
\special{pn 20}%
\special{sh 1}%
\special{ar 2400 1590 10 10 0  6.28318530717959E+0000}%
\special{sh 1}%
\special{ar 2400 1590 10 10 0  6.28318530717959E+0000}%
%
\special{pn 8}%
\special{ar 2400 1600 494 494  0.0128198 6.2831853}%
%
\special{pn 8}%
\special{ar 2400 1600 708 708  0.0000000 6.2831853}%
%
\special{pn 8}%
\special{ar 2400 1120 170 220  3.5102886 5.8175397}%
%
\special{pn 8}%
\special{ar 2780 900 250 310  1.2570314 2.7589330}%
%
\special{pn 8}%
\special{ar 2840 1330 190 130  4.7123890 6.2831853}%
\special{ar 2840 1330 190 130  0.0000000 0.7940186}%
%
\special{pn 8}%
\special{ar 3080 1590 180 200  2.2142974 4.1139701}%
%
\special{pn 8}%
\special{ar 2870 1860 140 160  5.6580178 6.2831853}%
\special{ar 2870 1860 140 160  0.0000000 1.8640993}%
%
\special{pn 8}%
\special{ar 2850 2180 350 160  3.0658756 4.6509550}%
%
\special{pn 8}%
\special{ar 2300 2130 210 160  0.4589992 3.1848989}%
%
\special{pn 8}%
\special{ar 1830 2200 260 380  4.7570686 6.0293236}%
%
\special{pn 8}%
\special{ar 1910 1560 200 260  1.8908885 4.5714850}%
%
\special{pn 8}%
\special{ar 1670 950 590 390  0.1680014 1.2070037}%
%
\special{pn 8}%
\special{pa 2990 840}%
\special{pa 2710 1090}%
\special{da 0.070}%
\special{sh 1}%
\special{pa 2710 1090}%
\special{pa 2774 1062}%
\special{pa 2750 1054}%
\special{pa 2746 1032}%
\special{pa 2710 1090}%
\special{fp}%
\put(30.5000,-8.6000){\makebox(0,0)[lb]{$\mathbb S^1(r_1^{-2})$}}%
%
\special{pn 8}%
\special{pa 2970 500}%
\special{pa 2510 960}%
\special{da 0.070}%
\special{sh 1}%
\special{pa 2510 960}%
\special{pa 2572 928}%
\special{pa 2548 922}%
\special{pa 2544 900}%
\special{pa 2510 960}%
\special{fp}%
\put(30.3000,-5.2000){\makebox(0,0)[lb]{$f_2(\mathbb S^1(r_2^{-2}))$}}%
%
\special{pn 8}%
\special{pa 3270 1890}%
\special{pa 3080 1830}%
\special{da 0.070}%
\special{sh 1}%
\special{pa 3080 1830}%
\special{pa 3138 1870}%
\special{pa 3132 1846}%
\special{pa 3150 1832}%
\special{pa 3080 1830}%
\special{fp}%
%
\special{pn 8}%
\special{pa 3200 2190}%
\special{pa 2800 1900}%
\special{da 0.070}%
\special{sh 1}%
\special{pa 2800 1900}%
\special{pa 2842 1956}%
\special{pa 2844 1932}%
\special{pa 2866 1924}%
\special{pa 2800 1900}%
\special{fp}%
\put(32.6000,-21.8000){\makebox(0,0)[lt]{$\mathbb S^1((r_1-\varepsilon)^{-2})$}}%
\put(33.1000,-18.4000){\makebox(0,0)[lt]{$\mathbb S^1((r_1+\varepsilon)^{-2})$}}%
\put(24.4000,-15.7000){\makebox(0,0)[lb]{$o$}}%
%
\special{pn 8}%
\special{pa 1470 1600}%
\special{pa 3490 1600}%
\special{fp}%
\special{sh 1}%
\special{pa 3490 1600}%
\special{pa 3424 1580}%
\special{pa 3438 1600}%
\special{pa 3424 1620}%
\special{pa 3490 1600}%
\special{fp}%
%
\special{pn 8}%
\special{pa 2400 2390}%
\special{pa 2400 630}%
\special{fp}%
\special{sh 1}%
\special{pa 2400 630}%
\special{pa 2380 698}%
\special{pa 2400 684}%
\special{pa 2420 698}%
\special{pa 2400 630}%
\special{fp}%
\put(8.1000,-25.6000){\makebox(0,0)[lt]{{\small Increase the number of waves of $f_2(\mathbb S^1(r_2^{-2}))$ as approach $\varepsilon$ to $0$.}}}%
\end{picture}%
\hspace{3truecm}}

\vspace{0.01truecm}

\centerline{{\bf Figure 2$\,\,:\,\,$ The first example showing 
$\breve d^k_{GH}([\mathbb S^1(r_1^{-2})],[\mathbb S^1(r_2^{-2})])=0$}}

\vspace{1truecm}

\centerline{
\unitlength 0.1in
\begin{picture}( 20.3000, 11.8000)( 15.0000,-25.2000)
%
\special{pn 8}%
\special{ar 2620 2050 416 96  0.0000000 6.2831853}%
%
\special{pn 8}%
\special{ar 2610 1770 350 70  4.7123890 6.2831853}%
\special{ar 2610 1770 350 70  0.0000000 1.5707963}%
%
\special{pn 8}%
\special{ar 2610 1800 420 100  3.7037794 4.7123890}%
%
\special{pn 8}%
\special{ar 2600 1810 430 110  1.5707963 3.5170486}%
%
\special{pn 8}%
\special{ar 2600 1750 430 170  3.4885095 6.2831853}%
\special{ar 2600 1750 430 170  0.0000000 1.5707963}%
\put(31.0000,-19.9000){\makebox(0,0)[lt]{$\widehat f_1(\mathbb S^1(r_1^{-2}))$}}%
\put(31.0000,-16.7000){\makebox(0,0)[lt]{$\widehat f_1(\mathbb S^1(r_2^{-2}))$}}%
%
\special{pn 8}%
\special{ar 2610 1710 410 130  1.5707963 3.3302319}%
%
\special{pn 8}%
\special{pa 1840 1920}%
\special{pa 2080 1820}%
\special{fp}%
\special{sh 1}%
\special{pa 2080 1820}%
\special{pa 2012 1828}%
\special{pa 2032 1842}%
\special{pa 2026 1864}%
\special{pa 2080 1820}%
\special{fp}%
%
\special{pn 8}%
\special{pa 1840 1920}%
\special{pa 2130 2010}%
\special{fp}%
\special{sh 1}%
\special{pa 2130 2010}%
\special{pa 2072 1972}%
\special{pa 2080 1994}%
\special{pa 2060 2010}%
\special{pa 2130 2010}%
\special{fp}%
\put(18.8000,-23.1000){\makebox(0,0)[lt]{{\small sufficiently close}}}%
%
\special{pn 8}%
\special{ar 1840 2140 340 220  1.5707963 4.7123890}%
\put(27.8000,-25.2000){\makebox(0,0)[lt]{$\,$}}%
%
\special{pn 8}%
\special{pa 3300 1540}%
\special{pa 3530 1540}%
\special{fp}%
%
\special{pn 8}%
\special{pa 3300 1540}%
\special{pa 3300 1340}%
\special{fp}%
\put(33.3000,-15.1000){\makebox(0,0)[lb]{$\mathbb E^3$}}%
\end{picture}%
\hspace{1truecm}}

\vspace{0.1truecm}

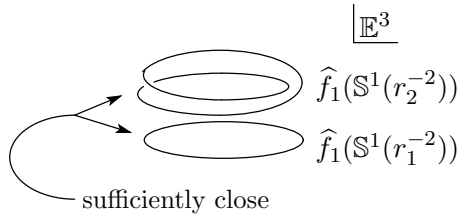
\centerline{{\bf Figure 3$\,\,:\,\,$ The second example showing 
$\breve d^k_{GH}([\mathbb S^1(r_1^{-2})],[\mathbb S^1(r_2^{-2})])=0$}}

\vspace{0.5truecm}

\noindent
We shall give third example showing that $\breve d^k_{GH}$ is not a distance function.  
Let $\mathbb R^{n+2}_1$ be the $(n+2)$-dimensional Lorentzian space and $g_{\mathbb L}$ the Lorenzian metric of 
$\mathbb R^{n+2}_1$, that is, $g_{\mathbb L}=-dx_1^2+dx_2^2+\cdots+dx_{n+2}^2$.  
Put 
$$\mathbb H^{n+1}(-\widetilde r^{-2}):=\{(x_1,\cdots,x_{n+2})\in\mathbb R^{n+2}_1\,|\,
-x_1^2+x_2^2+\cdots+x_{n+2}^2=-\widetilde r^2\}\quad\,\,(\widetilde r>0).$$
Denote by $\iota^H_{\widetilde r}$ the inclusion map of $\mathbb H^{n+1}(-\widetilde r^{-2})$ into 
$\mathbb R^{n+2}_1$ and $g^H_{\widetilde r}$ the induced metric $(\iota^H_{\widetilde r})^{\ast}g_{\mathbb L}$.  
The space $(\mathbb H^{n+1}(-{\widetilde r}^{-2}),g^H_{\widetilde r})$ is the $(n+1)$-dimensional hyperbolic 
space of constant curvature $-{\widetilde r}^{-2}$.  The sphere $(\mathbb S^n(r^{-2}),g^S_r)$ is isometrically 
embedded into $(\mathbb H^{n+1}(-\widetilde r^{-2}),g^H_{\widetilde r})$ by the following $C^{\infty}$-embedding:
$$f_{r,\widetilde r}(x_1,\cdots,x_{n+1})=(\sqrt{\widetilde r^2+r^2},x_1,\cdots,x_{n+1})\quad\,\,
((x_1,\cdots,x_{n+1})\in \mathbb S^n(r^{-2})).$$
Take distinct positive constants $r_1$ and $r_2$ ($r_1<r_2$).  
We shall calculate 
$$d_{H,\mathbb H^{n+1}(-\widetilde r^{-2})}(f_{r_1,\widetilde r}(\mathbb S^n(r_1^{-2})),
f_{r_2,\widetilde r}(\mathbb S^n(r_2^{-2}))),$$
which is equal to 
\begin{align*}
&d_{g^H_{\widetilde r}}(f_{r_1,\widetilde r}(r_1,0,\cdots,0),f_{r_2,\widetilde r}(r_2,0,\cdots,0))\\
(=&d_{g^H_{\widetilde r}}((\sqrt{\widetilde r^2+r_1^2},r_1,0,\cdots,0),
(\sqrt{\widetilde r^2+r_2^2},r_2,0,\cdots,0))).
\end{align*}
The shortest geodesic $\gamma_{\widetilde r}$ (in $\mathbb H^{n+1}(-\widetilde r^{-2})$) connecting 
$(\sqrt{\widetilde r^2+r_1^2},r_1,0,\cdots,0)$ and \newline
$(\sqrt{\widetilde r^2+r_2^2},r_2,0,\cdots,0)$ is given by 
$$\gamma_{\widetilde r}(t):=(\widetilde r\,\cosh\,t,\widetilde r\,\sinh\,t,0,\cdots,0)\quad\,\,
(\sinh^{-1}\left(\frac{r_1}{\widetilde r}\right)\leq t\leq\sinh^{-1}\left(\frac{r_2}{\widetilde r}\right)),$$
where $\sinh^{-1}$ denotes the inverse function of $\sinh|_{[0,\infty)}$.  
For the simplicity, put $\displaystyle{a(r_i):=\sinh^{-1}\left(\frac{r_i}{\widetilde r}\right)}$ ($i=1,2$).  
The length $L(\gamma_{\widetilde r})$ of $\gamma$ is given by 
\begin{align*}
L(\gamma_{\widetilde r})&=\int_{a(r_1)}^{a(r_2)}\|\gamma'_{\widetilde r}(t)\|\,dt
=\int_{a(r_1)}^{a(r_2)}\sqrt{|-\widetilde r^2|}\,dt\\
&=\widetilde r\cdot(a(r_2)-a(r_1))=\widetilde r\cdot
\left(\sinh^{-1}\left(\frac{r_2}{\widetilde r}\right)-\sinh^{-1}\left(\frac{r_1}{\widetilde r}\right)\right).
\end{align*}
For the simplicity, denote by $F_{r_1,r_2}(\widetilde r)$ the right-hand side of this relation.  
Then we have 
$$d_{g^H_{\widetilde r}}(f_{r_1,\widetilde r}(r_1,0,\cdots,0),f_{r_1,\widetilde r}(r_2,0,\cdots,0))
=F_{r_1,r_2}(\widetilde r),$$
that is, 
$$d_{H,\mathbb H^{n+1}(-\widetilde r^{-2})}(f_{r_1,\widetilde r}(\mathbb S^n(r_1^{-2})),
f_{r_2,\widetilde r}(\mathbb S^n(r_2^{-2})))=F_{r_1,r_2}(\widetilde r).\leqno{(2.1)}$$
By using L'H$\hat{\rm o}$pital's theorem, we have 
$$\begin{array}{l}
\displaystyle{\lim_{\widetilde r\to +0}F_{r_1,r_2}(\widetilde r)
=\lim_{\widetilde r\to +0}\left(\left(1+\left(\frac{r_2}{\widetilde r}\right)^2\right)^{-\frac{1}{2}}
-\left(1+\left(\frac{r_1}{\widetilde r}\right)^2\right)^{-\frac{1}{2}}\right)}\\
\hspace{2.375truecm}\displaystyle{=\lim_{\widetilde r\to +0}
\frac{r_2\sqrt{\widetilde r^2+r_1^2}-r_1\sqrt{\widetilde r^2+r_2^2}}
{\sqrt{\widetilde r^2+r_1^2}\cdot\sqrt{\widetilde r^2+r_2^2}}=0.}
\end{array}\leqno{(2.2)}$$
Hence we obtain 
$$\lim_{\widetilde r\to +0}d_{H,\mathbb H^{n+1}(-\widetilde r^{-2})}(f_{r_1,\widetilde r}(\mathbb S^n(r_1^{-2})),
f_{r_2,\widetilde r}(\mathbb S^n(r_2^{-2})))=0,$$
that is, 
$$\breve d^{\infty}_{GH}([\mathbb S^n(r_1^{-2})],[\mathbb S^n(r_2^{-2})])=0.$$
In more general, we can give the following counter-examples.  Let $M$ be a compact submanifold in 
$\mathbb S^n(1)$ embedded by a $C^{k+1}$-embedding $f$.  Then, since $(M,f^{\ast}g^S_1)$ and 
$(M,r\cdot f^{\ast}g^S_1)$ ($r>0,r\not=1$) 
are compact $C^k$-Riemannian submanifolds in $\mathbb S^n(1)$ and $\mathbb S^n(r^{-2})$, respectively, 
it is shown that 
$$d_{H,\mathbb H^{n+1}(-\widetilde r^{-2})}((M,f^{\ast}g^S_1),(M,r\cdot f^{\ast}g^S_1))
=F_{r_1,r_2}(\widetilde r)$$
and hence 
$$\breve d^k_{GH}([(M,f^{\ast}g^S_1)],[(M,r\cdot f^{\ast}g^S_1)])=0.$$
On the other hand, it is clear that $[(M,f^{\ast}g^S_1)]\not=[(M,r\cdot f^{\ast}g^S_1)]$.  
Denote by $D_{r_1,r_2}$ the domain of $\mathbb H^{n+1}(-\widetilde r^{-2})$ sorrounded by 
$f_{r_1\widetilde r}(\mathbb S^n(r_1^{-2}))$ and $f_{r_2\widetilde r}(\mathbb S^n(r_2^{-2}))$.  
The Riemannian manifold $(D_{r_1,r_2},g^H_{\widetilde r}|_{D_{r_1,r_2}})$ cannot be isometrically embedded into 
$\mathbb E^m$ for any $m\in\mathbb N$ (see Figure 5).  We consider that this fact arises 
$$\lim_{\widetilde r\to +0}d_{H,\mathbb H^{n+1}(-\widetilde r^{-2})}(f_{r_1,\widetilde r}(\mathbb S^n(r_1^{-2})),
f_{r_2,\widetilde r}(\mathbb S^n(r_2^{-2})))=0.$$

\vspace{0.3truecm}

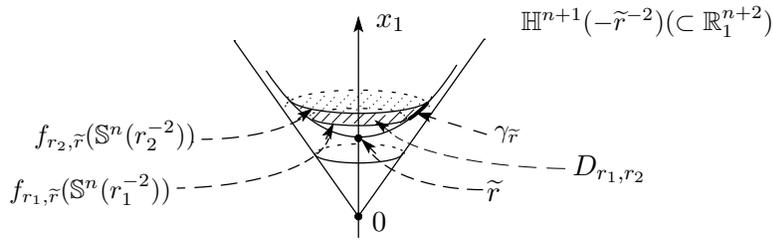
\centerline{
\unitlength 0.1in
\begin{picture}( 69.6000, 15.9000)(-29.5000,-24.5000)
%
\special{pn 8}%
\special{ar 2884 1472 246 38  3.1415927 3.2263983}%
\special{ar 2884 1472 246 38  3.4808153 3.5656209}%
\special{ar 2884 1472 246 38  3.8200379 3.9048435}%
\special{ar 2884 1472 246 38  4.1592605 4.2440662}%
\special{ar 2884 1472 246 38  4.4984831 4.5832888}%
\special{ar 2884 1472 246 38  4.8377057 4.9225114}%
\special{ar 2884 1472 246 38  5.1769283 5.2617340}%
\special{ar 2884 1472 246 38  5.5161510 5.6009566}%
\special{ar 2884 1472 246 38  5.8553736 5.9401792}%
\special{ar 2884 1472 246 38  6.1945962 6.2794018}%
%
\special{pn 8}%
\special{ar 2896 1374 376 86  6.2831853 6.2831853}%
\special{ar 2896 1374 376 86  0.0000000 3.1415927}%
%
\special{pn 8}%
\special{ar 2884 1424 376 86  3.1415927 3.1935407}%
\special{ar 2884 1424 376 86  3.3493849 3.4013329}%
\special{ar 2884 1424 376 86  3.5571771 3.6091251}%
\special{ar 2884 1424 376 86  3.7649693 3.8169173}%
\special{ar 2884 1424 376 86  3.9727615 4.0247095}%
\special{ar 2884 1424 376 86  4.1805537 4.2325017}%
\special{ar 2884 1424 376 86  4.3883459 4.4402940}%
\special{ar 2884 1424 376 86  4.5961381 4.6480862}%
\special{ar 2884 1424 376 86  4.8039303 4.8558784}%
\special{ar 2884 1424 376 86  5.0117225 5.0636706}%
\special{ar 2884 1424 376 86  5.2195147 5.2714628}%
\special{ar 2884 1424 376 86  5.4273069 5.4792550}%
\special{ar 2884 1424 376 86  5.6350991 5.6870472}%
\special{ar 2884 1424 376 86  5.8428914 5.8948394}%
\special{ar 2884 1424 376 86  6.0506836 6.1026316}%
\special{ar 2884 1424 376 86  6.2584758 6.2831853}%
%
\special{pn 8}%
\special{ar 2896 1486 246 38  6.2831853 6.2831853}%
\special{ar 2896 1486 246 38  0.0000000 3.1415927}%
%
\special{pn 8}%
\special{pa 2896 1990}%
\special{pa 3548 1072}%
\special{fp}%
\special{pa 3548 1072}%
\special{pa 3548 1072}%
\special{fp}%
%
\special{pn 8}%
\special{pa 2890 2000}%
\special{pa 2240 1080}%
\special{fp}%
\special{pa 2240 1080}%
\special{pa 2240 1080}%
\special{fp}%
%
\special{pn 8}%
\special{ar 2896 -1412 1046 2996  1.0642157 2.0561297}%
%
\special{pn 8}%
\special{ar 2064 1374 612 222  0.5487308 0.6929615}%
\special{ar 2064 1374 612 222  0.7795000 0.9237308}%
\special{ar 2064 1374 612 222  1.0102692 1.1545000}%
\special{ar 2064 1374 612 222  1.2410385 1.3852692}%
\special{ar 2064 1374 612 222  1.4718077 1.5594332}%
%
\special{pn 8}%
\special{pa 2574 1490}%
\special{pa 2644 1440}%
\special{da 0.070}%
\special{sh 1}%
\special{pa 2644 1440}%
\special{pa 2578 1462}%
\special{pa 2600 1472}%
\special{pa 2600 1496}%
\special{pa 2644 1440}%
\special{fp}%
%
\special{pn 8}%
\special{ar 1906 1556 834 296  0.1201826 0.2264713}%
\special{ar 1906 1556 834 296  0.2902446 0.3965333}%
\special{ar 1906 1556 834 296  0.4603066 0.5665953}%
\special{ar 1906 1556 834 296  0.6303686 0.7366573}%
\special{ar 1906 1556 834 296  0.8004306 0.9067193}%
\special{ar 1906 1556 834 296  0.9704926 1.0767813}%
\special{ar 1906 1556 834 296  1.1405546 1.2468433}%
\special{ar 1906 1556 834 296  1.3106166 1.4169053}%
\special{ar 1906 1556 834 296  1.4806786 1.5286310}%
%
\special{pn 8}%
\special{pa 2720 1594}%
\special{pa 2756 1520}%
\special{da 0.070}%
\special{sh 1}%
\special{pa 2756 1520}%
\special{pa 2708 1572}%
\special{pa 2732 1568}%
\special{pa 2744 1590}%
\special{pa 2756 1520}%
\special{fp}%
\put(20.4000,-16.6000){\makebox(0,0)[rb]{{\small $f_{r_2,\widetilde r}(\mathbb S^n(r_2^{-2}))$}}}%
\put(37.4000,-10.6000){\makebox(0,0)[lb]{{\small $\mathbb H^{n+1}(-\widetilde r^{-2})(\subset\mathbb R^{n+2}_1)$}}}%
\put(7.1000,-22.6000){\makebox(0,0)[lt]{{\small $d_{G,\mathbb H^{n+1}(-\widetilde r^{-2})}(f_{r_1,\widetilde r}(\mathbb S^n(r_1^{-2})),\,f_{r_2,\widetilde r}(\mathbb S^n(r_2^{-2})))=\widetilde r\cdot(a(r_2)-a(r_1))\,\,\to\,\,0\quad(\widetilde r\to +0)$}}}%
\put(36.0000,-15.3000){\makebox(0,0)[lt]{{\small $\gamma_{\widetilde r}$}}}%
\put(23.4000,-24.5000){\makebox(0,0)[lt]{$\,$}}%
\put(19.1000,-17.8000){\makebox(0,0)[rt]{{\small $f_{r_1,\widetilde r}(\mathbb S^n(r_1^{-2}))$}}}%
%
\special{pn 8}%
\special{ar 2890 1670 230 50  6.2831853 6.2831853}%
\special{ar 2890 1670 230 50  0.0000000 3.1415927}%
%
\special{pn 8}%
\special{ar 2900 1670 230 50  3.1415927 3.2273069}%
\special{ar 2900 1670 230 50  3.4844498 3.5701641}%
\special{ar 2900 1670 230 50  3.8273069 3.9130212}%
\special{ar 2900 1670 230 50  4.1701641 4.2558784}%
\special{ar 2900 1670 230 50  4.5130212 4.5987355}%
\special{ar 2900 1670 230 50  4.8558784 4.9415927}%
\special{ar 2900 1670 230 50  5.1987355 5.2844498}%
\special{ar 2900 1670 230 50  5.5415927 5.6273069}%
\special{ar 2900 1670 230 50  5.8844498 5.9701641}%
\special{ar 2900 1670 230 50  6.2273069 6.2831853}%
%
\special{pn 8}%
\special{pa 2890 2110}%
\special{pa 2890 960}%
\special{fp}%
\special{sh 1}%
\special{pa 2890 960}%
\special{pa 2870 1028}%
\special{pa 2890 1014}%
\special{pa 2910 1028}%
\special{pa 2890 960}%
\special{fp}%
\put(29.8000,-10.3000){\makebox(0,0)[lb]{$x_1$}}%
%
\special{pn 20}%
\special{sh 1}%
\special{ar 2890 2000 10 10 0  6.28318530717959E+0000}%
\special{sh 1}%
\special{ar 2890 2000 10 10 0  6.28318530717959E+0000}%
%
\special{pn 20}%
\special{sh 1}%
\special{ar 2890 1590 10 10 0  6.28318530717959E+0000}%
\special{sh 1}%
\special{ar 2890 1590 10 10 0  6.28318530717959E+0000}%
\put(29.6000,-19.8000){\makebox(0,0)[lt]{$0$}}%
\put(35.6000,-18.1000){\makebox(0,0)[lt]{$\widetilde r$}}%
%
\special{pn 8}%
\special{pa 2020 1610}%
\special{pa 2020 1610}%
\special{dt 0.045}%
%
\special{pn 8}%
\special{ar 3540 1550 630 310  1.5707963 1.6984559}%
\special{ar 3540 1550 630 310  1.7750516 1.9027112}%
\special{ar 3540 1550 630 310  1.9793070 2.1069665}%
\special{ar 3540 1550 630 310  2.1835623 2.3112219}%
\special{ar 3540 1550 630 310  2.3878176 2.5154772}%
\special{ar 3540 1550 630 310  2.5920729 2.7197325}%
%
\special{pn 8}%
\special{pa 2980 1670}%
\special{pa 2910 1600}%
\special{da 0.070}%
\special{sh 1}%
\special{pa 2910 1600}%
\special{pa 2944 1662}%
\special{pa 2948 1638}%
\special{pa 2972 1634}%
\special{pa 2910 1600}%
\special{fp}%
%
\special{pn 8}%
\special{pa 3560 1570}%
\special{pa 3200 1460}%
\special{da 0.070}%
\special{sh 1}%
\special{pa 3200 1460}%
\special{pa 3258 1500}%
\special{pa 3252 1476}%
\special{pa 3270 1460}%
\special{pa 3200 1460}%
\special{fp}%
%
\special{pn 20}%
\special{ar 2890 1040 500 530  0.7704739 1.0043899}%
%
\special{pn 4}%
\special{pa 2586 1436}%
\special{pa 2596 1426}%
\special{fp}%
\special{pa 2614 1468}%
\special{pa 2638 1442}%
\special{fp}%
\special{pa 2646 1496}%
\special{pa 2688 1452}%
\special{fp}%
\special{pa 2700 1502}%
\special{pa 2740 1460}%
\special{fp}%
\special{pa 2752 1510}%
\special{pa 2800 1460}%
\special{fp}%
\special{pa 2810 1510}%
\special{pa 2860 1460}%
\special{fp}%
\special{pa 2868 1512}%
\special{pa 2920 1460}%
\special{fp}%
\special{pa 2920 1520}%
\special{pa 2978 1462}%
\special{fp}%
\special{pa 2980 1520}%
\special{pa 3030 1470}%
\special{fp}%
\special{pa 3050 1510}%
\special{pa 3110 1450}%
\special{fp}%
\special{pa 3110 1510}%
\special{pa 3180 1440}%
\special{fp}%
%
\special{pn 4}%
\special{pa 2600 1420}%
\special{pa 2650 1370}%
\special{dt 0.027}%
\special{pa 2652 1428}%
\special{pa 2714 1368}%
\special{dt 0.027}%
\special{pa 2706 1436}%
\special{pa 2774 1366}%
\special{dt 0.027}%
\special{pa 2758 1442}%
\special{pa 2840 1362}%
\special{dt 0.027}%
\special{pa 2810 1450}%
\special{pa 2910 1350}%
\special{dt 0.027}%
\special{pa 2880 1442}%
\special{pa 2966 1354}%
\special{dt 0.027}%
\special{pa 2948 1432}%
\special{pa 3022 1360}%
\special{dt 0.027}%
\special{pa 3004 1438}%
\special{pa 3076 1366}%
\special{dt 0.027}%
\special{pa 3060 1440}%
\special{pa 3130 1370}%
\special{dt 0.027}%
\special{pa 3126 1436}%
\special{pa 3176 1384}%
\special{dt 0.027}%
\special{pa 3190 1430}%
\special{pa 3222 1398}%
\special{dt 0.027}%
%
\special{pn 4}%
\special{ar 4040 1480 1030 270  1.6378997 1.7302074}%
\special{ar 4040 1480 1030 270  1.7855920 1.8778997}%
\special{ar 4040 1480 1030 270  1.9332843 2.0255920}%
\special{ar 4040 1480 1030 270  2.0809766 2.1732843}%
\special{ar 4040 1480 1030 270  2.2286690 2.3209766}%
\special{ar 4040 1480 1030 270  2.3763613 2.4686690}%
\special{ar 4040 1480 1030 270  2.5240536 2.6163613}%
\special{ar 4040 1480 1030 270  2.6717459 2.7640536}%
\special{ar 4040 1480 1030 270  2.8194382 2.8374110}%
%
\special{pn 4}%
\special{pa 3060 1560}%
\special{pa 2980 1500}%
\special{da 0.070}%
\special{sh 1}%
\special{pa 2980 1500}%
\special{pa 3022 1556}%
\special{pa 3024 1532}%
\special{pa 3046 1524}%
\special{pa 2980 1500}%
\special{fp}%
\put(40.1000,-16.8000){\makebox(0,0)[lt]{$D_{r_1,r_2}$}}%
\end{picture}%
\hspace{14truecm}}

\vspace{0.3truecm}

\centerline{{\bf Figure 4$\,\,:\,\,$ The third example showing $\breve d^{\infty}_{GH}([\mathbb S^n(r_1^{-2})],
[\mathbb S^n(r_2^{-2})])=0$}}

\vspace{1truecm}

\centerline{
\unitlength 0.1in
\begin{picture}( 84.8000, 17.4000)(-32.8000,-23.3000)
\put(29.9000,-23.3000){\makebox(0,0)[lt]{$\,$}}%
%
\special{pn 8}%
\special{ar 2350 1302 210 30  3.1415927 3.2415927}%
\special{ar 2350 1302 210 30  3.5415927 3.6415927}%
\special{ar 2350 1302 210 30  3.9415927 4.0415927}%
\special{ar 2350 1302 210 30  4.3415927 4.4415927}%
\special{ar 2350 1302 210 30  4.7415927 4.8415927}%
\special{ar 2350 1302 210 30  5.1415927 5.2415927}%
\special{ar 2350 1302 210 30  5.5415927 5.6415927}%
\special{ar 2350 1302 210 30  5.9415927 6.0415927}%
%
\special{pn 8}%
\special{ar 2340 1222 320 70  6.2831853 6.2831853}%
\special{ar 2340 1222 320 70  0.0000000 3.1415927}%
%
\special{pn 8}%
\special{ar 2350 1262 320 70  3.1415927 3.2031311}%
\special{ar 2350 1262 320 70  3.3877465 3.4492850}%
\special{ar 2350 1262 320 70  3.6339003 3.6954388}%
\special{ar 2350 1262 320 70  3.8800542 3.9415927}%
\special{ar 2350 1262 320 70  4.1262080 4.1877465}%
\special{ar 2350 1262 320 70  4.3723619 4.4339003}%
\special{ar 2350 1262 320 70  4.6185157 4.6800542}%
\special{ar 2350 1262 320 70  4.8646696 4.9262080}%
\special{ar 2350 1262 320 70  5.1108234 5.1723619}%
\special{ar 2350 1262 320 70  5.3569773 5.4185157}%
\special{ar 2350 1262 320 70  5.6031311 5.6646696}%
\special{ar 2350 1262 320 70  5.8492850 5.9108234}%
\special{ar 2350 1262 320 70  6.0954388 6.1569773}%
%
\special{pn 8}%
\special{ar 2340 1312 210 30  6.2831853 6.2831853}%
\special{ar 2340 1312 210 30  0.0000000 3.1415927}%
%
\special{pn 8}%
\special{pa 2340 1722}%
\special{pa 1790 970}%
\special{dt 0.045}%
\special{pa 1790 970}%
\special{pa 1790 970}%
\special{dt 0.045}%
%
\special{pn 8}%
\special{ar 2340 -1038 890 2430  1.0863184 2.0775147}%
%
\special{pn 8}%
\special{ar 1620 1222 520 180  0.5468445 0.7182731}%
\special{ar 1620 1222 520 180  0.8211302 0.9925588}%
\special{ar 1620 1222 520 180  1.0954159 1.2668445}%
\special{ar 1620 1222 520 180  1.3697016 1.5411302}%
%
\special{pn 8}%
\special{pa 2050 1322}%
\special{pa 2110 1282}%
\special{da 0.070}%
\special{sh 1}%
\special{pa 2110 1282}%
\special{pa 2044 1302}%
\special{pa 2066 1312}%
\special{pa 2066 1336}%
\special{pa 2110 1282}%
\special{fp}%
%
\special{pn 8}%
\special{ar 1590 1342 710 240  0.1198313 0.2461471}%
\special{ar 1590 1342 710 240  0.3219365 0.4482523}%
\special{ar 1590 1342 710 240  0.5240418 0.6503576}%
\special{ar 1590 1342 710 240  0.7261471 0.8524629}%
\special{ar 1590 1342 710 240  0.9282523 1.0545681}%
\special{ar 1590 1342 710 240  1.1303576 1.2566734}%
\special{ar 1590 1342 710 240  1.3324629 1.4587786}%
%
\special{pn 8}%
\special{pa 2270 1402}%
\special{pa 2300 1342}%
\special{da 0.070}%
\special{sh 1}%
\special{pa 2300 1342}%
\special{pa 2252 1392}%
\special{pa 2276 1390}%
\special{pa 2288 1410}%
\special{pa 2300 1342}%
\special{fp}%
\put(15.8000,-14.7000){\makebox(0,0)[rb]{{\small $f_{r_2,\widetilde r}(\mathbb S^n(r_2^{-2}))$}}}%
\put(27.8000,-7.6000){\makebox(0,0)[lb]{{\small $\mathbb H^{n+1}(-\widetilde r^{-2})\,(\subset\mathbb R^{n+2}_1)$}}}%
%
\special{pn 8}%
\special{pa 2350 1722}%
\special{pa 2900 970}%
\special{dt 0.045}%
\special{pa 2900 970}%
\special{pa 2900 970}%
\special{dt 0.045}%
%
\special{pn 8}%
\special{ar 2350 1512 160 40  6.2831853 6.4031853}%
\special{ar 2350 1512 160 40  6.7631853 6.8831853}%
\special{ar 2350 1512 160 40  7.2431853 7.3631853}%
\special{ar 2350 1512 160 40  7.7231853 7.8431853}%
\special{ar 2350 1512 160 40  8.2031853 8.3231853}%
\special{ar 2350 1512 160 40  8.6831853 8.8031853}%
\special{ar 2350 1512 160 40  9.1631853 9.2831853}%
%
\special{pn 8}%
\special{ar 2350 1512 160 40  3.1415927 3.2615927}%
\special{ar 2350 1512 160 40  3.6215927 3.7415927}%
\special{ar 2350 1512 160 40  4.1015927 4.2215927}%
\special{ar 2350 1512 160 40  4.5815927 4.7015927}%
\special{ar 2350 1512 160 40  5.0615927 5.1815927}%
\special{ar 2350 1512 160 40  5.5415927 5.6615927}%
\special{ar 2350 1512 160 40  6.0215927 6.1415927}%
%
\special{pn 8}%
\special{ar 2060 830 400 790  0.3740565 0.7914587}%
%
\special{pn 8}%
\special{ar 1830 830 860 710  0.4354045 0.9317683}%
%
\special{pn 8}%
\special{ar 2360 1400 70 50  5.3326385 5.3783582}%
%
\special{pn 8}%
\special{pa 2360 1400}%
\special{pa 2356 1396}%
\special{fp}%
\special{pa 2396 1376}%
\special{pa 2370 1350}%
\special{fp}%
\special{pa 2430 1350}%
\special{pa 2394 1314}%
\special{fp}%
\special{pa 2456 1316}%
\special{pa 2410 1270}%
\special{fp}%
\special{pa 2482 1282}%
\special{pa 2418 1218}%
\special{fp}%
\special{pa 2510 1250}%
\special{pa 2432 1172}%
\special{fp}%
\special{pa 2540 1220}%
\special{pa 2452 1132}%
\special{fp}%
\special{pa 2570 1190}%
\special{pa 2520 1140}%
\special{fp}%
\special{pa 2596 1156}%
\special{pa 2590 1150}%
\special{fp}%
%
\special{pn 4}%
\special{ar 3930 1710 250 430  5.4312190 6.2831853}%
%
\special{pn 4}%
\special{ar 4420 1710 250 430  3.1415927 3.9935590}%
%
\special{pn 4}%
\special{ar 4180 1380 80 30  0.0000000 6.2831853}%
\put(51.7000,-13.5000){\makebox(0,0)[lb]{$\mathbb E^{n+2}$}}%
\put(38.2000,-15.0000){\makebox(0,0)[rt]{{\small isometric embedding}}}%
%
\special{pn 8}%
\special{pa 4170 1630}%
\special{pa 4158 1618}%
\special{fp}%
\special{pa 4182 1582}%
\special{pa 4156 1556}%
\special{fp}%
\special{pa 4194 1534}%
\special{pa 4130 1470}%
\special{fp}%
\special{pa 4206 1486}%
\special{pa 4098 1378}%
\special{fp}%
\special{pa 4216 1436}%
\special{pa 4138 1358}%
\special{fp}%
\special{pa 4236 1396}%
\special{pa 4196 1356}%
\special{fp}%
%
\special{pn 4}%
\special{pa 2630 1200}%
\special{pa 3950 1500}%
\special{fp}%
\special{sh 1}%
\special{pa 3950 1500}%
\special{pa 3890 1466}%
\special{pa 3898 1488}%
\special{pa 3882 1506}%
\special{pa 3950 1500}%
\special{fp}%
%
\special{pn 4}%
\special{pa 3610 1880}%
\special{pa 5200 1880}%
\special{fp}%
\special{sh 1}%
\special{pa 5200 1880}%
\special{pa 5134 1860}%
\special{pa 5148 1880}%
\special{pa 5134 1900}%
\special{pa 5200 1880}%
\special{fp}%
%
\special{pn 4}%
\special{pa 4090 2160}%
\special{pa 5070 1550}%
\special{fp}%
\special{sh 1}%
\special{pa 5070 1550}%
\special{pa 5004 1568}%
\special{pa 5026 1578}%
\special{pa 5024 1602}%
\special{pa 5070 1550}%
\special{fp}%
%
\special{pn 4}%
\special{pa 4540 2200}%
\special{pa 4540 1260}%
\special{fp}%
\special{sh 1}%
\special{pa 4540 1260}%
\special{pa 4520 1328}%
\special{pa 4540 1314}%
\special{pa 4560 1328}%
\special{pa 4540 1260}%
\special{fp}%
\put(15.8000,-15.3000){\makebox(0,0)[rt]{{\small $f_{r_1,\widetilde r}(\mathbb S^n(r_1^{-2}))$}}}%
%
\special{pn 8}%
\special{ar 2690 790 340 370  0.0518055 0.2208195}%
\special{ar 2690 790 340 370  0.3222280 0.4912421}%
\special{ar 2690 790 340 370  0.5926505 0.7616646}%
\special{ar 2690 790 340 370  0.8630731 1.0320871}%
\special{ar 2690 790 340 370  1.1334956 1.3025097}%
\special{ar 2690 790 340 370  1.4039181 1.4305932}%
%
\special{pn 8}%
\special{pa 2750 1160}%
\special{pa 2670 1160}%
\special{da 0.070}%
\special{sh 1}%
\special{pa 2670 1160}%
\special{pa 2738 1180}%
\special{pa 2724 1160}%
\special{pa 2738 1140}%
\special{pa 2670 1160}%
\special{fp}%
\end{picture}%
\hspace{11truecm}}

\vspace{0.2truecm}

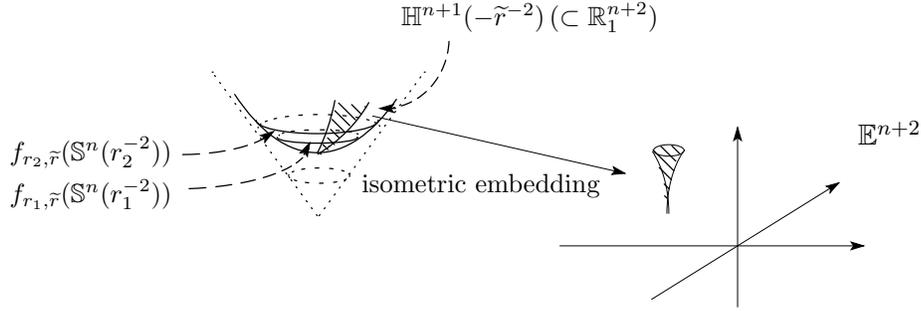
\centerline{{\bf Figure 5$\,\,:\,\,$ A domain in $\mathbb H^{n+1}(-\widetilde r^{-2})$ isometrically embedded into 
$\mathbb E^{n+2}$}}

\section{Gromov-Hausdorff-like distance function} 
By refering three examples in the previous section, we shall define a Gromov-Hausdorff-like distance function 
over $\mathcal{RM}_c^k$.  We use the notations in the previous section.  
The first example in the previous section indicates that $d_{H,(\widetilde M,d_{\widetilde g})}$ in the 
definitoin of $\breve d^k_{GH}$ should be replaced by a distance function including informations of the $i$-th 
derivatives ($1\leq i\leq k+1$) of the isometric embeddings becuase 
${\rm Emb}_I^{k+1}((M_i,g_i),(\widetilde M,\widetilde g))$ is a very wider class than 
${\rm Emb}_{d.p.}((M_i,d_{g_i}),(\widetilde M,d_{\widetilde g}))$.  
Also, the third example indicates that, in the definition of $\breve d^k_{GH}$, the range which 
$(\widetilde M,\widetilde g)$ moves should be restricted to the class of $C^k$-Riemannian manifolds 
isometrically embedded into a Euclidian space.  
On the basis of these reasons, we shall define another function over $\mathcal{RM}^k_c\times\mathcal{RM}^k_c$.  
We shall prepare some notions to state the definition.  
Let $\pi_1:TM\to M$ be the tangent bundle of $M$, $\pi_2:T(TM)\to TM$ be the tangent bundle of the manifold 
$TM$.  Denote by $T^2M$ the manifold $T(TM)$.  
Let $\pi_3:T(T^2M)\to T^2M$ be the tangent bundle of the manifold $T^2M$ and denote by $T^3M$ the manifold 
$T(T^2M)$.  In the sequel, we define $T^lM$ and $\pi_l$ ($l=4,5,\cdots$) inductively.  
Let $f$ be a $C^{k+1}$-map from a Riemannian manifold $(M,g)$ to another Riemannian manifold 
$(\widetilde M,\widetilde g)$.  
The differential $df:TM\to T\widetilde M$ of $f$ is defined by $df|_{T_xM}=df_x\,\,(x\in M)$ and 
the differential $d^2f:=d(df):T^2M\to T^2\widetilde M$ is defined similarly.  
In the sequel, $d^lf:T^lM\to T^l\widetilde M$ ($l=3,\cdots,k$) are defined inductively.  
The Sasaki metirc $g_S^1$ of $TM$ with respect to $g$ is defined by 
$$\begin{array}{r}
\displaystyle{(g_S^1)_{v_1}(v_2,v_3):=g_{\pi_1(v_1)}((d\pi_1)_{v_1}(v_2),(d\pi_1)_{v_1}(v_3))
+g_{\pi_1(v_1)}((v_2)_{\mathcal V},(v_3)_{\mathcal V})}\\
\displaystyle{\quad(v_1\in TM,\,\,v_2,v_3\in T_{v_1}(TM)),}
\end{array}$$
where $(v_i)_{\mathcal V}$ denotes the vertical component of $v_i$ with respect to 
$T_{v_1}(TM)=\mathcal V_{v_1}\oplus\mathcal H_{v_1}$ 
($\mathcal V:$ the vertical distribution, $\mathcal H:$ the horizontal distribution 
associated to the Riemannian connection of $g$).  Here we note that 
$\mathcal V_{v_1}(=T_{v_1}(\pi_1^{-1}(\pi_1(v_1))))$ is identified with 
$T_{v_1}(T_{\pi_1(v_1)}M)(=T_{\pi_1(v_1)}M)$.  
Similarly, the Sasaki metirc $g_S^2$ of $T^2M$ with respect to $g_S^1$ is defined.  
In the sequel, the Sasaki metirc $g_S^l$ of $T^lM$ with respect to $g_S^{l-1}$ ($l=3,4,\cdots,k$) are defined 
inductively.  Similarly, $T^l\widetilde M$ and $\widetilde g_S^l$ are defined for 
$(\widetilde M,\widetilde g)$.  Let $S^lM$ be the unit tangent bundle of the Riemannian manifold 
$(T^lM,g_S^{l-1})$ (i.e., $S^lM:=\{v\in T^lM\,|\,g_S^{l-1}(v,v)=1\}$).  
Denote by $m(M_i,g_i)$ the minimum of natural numbers $l$'s such that $(M_i,g_i)$ is isometrically 
embedded into $\mathbb E^l$.  For simplicity, set $m_i:=m(M_i,g_i)$.  
Here we note that, since $M_i$ is compact, the existence of such a minimum number is assured by 
the Nash's isometric embedding theorem (\cite{N1},\cite{N2}).   
On the basis of the above reasons and the consideration of the second example in the previous section, 
we define a function $d_{GH}^k:\mathcal{RM}^k_c\times\mathcal{RM}^k_c\to\mathbb R$ by 
{\small{
\begin{align*}
&d_{GH}^k([(M_1,g_1)],[(M_2,g_2)])\\
:=&\inf\{\sum_{j=0}^{k+1}\,d_{H,(T^j\mathbb E^m,d_{(g_{\mathbb E})_S^j})}
(d^j(\iota_1\circ f_1)(S^jM_1),d^j(\iota_2\circ f_2\circ\psi_2)(S^jM_2))\\
&\hspace{1.15truecm}\,|\,f_i\in{\rm Emb}_I^j((M_i,g_i),\mathbb E^{m_i}),\,\,\,\,
\iota_i\in{\rm Emb}^{t.g.}_I(\mathbb E^{m_i},\mathbb E^m)\,\,\,\,(i=1,2)),\\
&\hspace{1.35truecm}\psi_2\in{\rm Diff}^j(M_2)\},
\end{align*}}
where $m$ is a any natural number with $m\geq\max\{m_1,m_2\}$, ${\rm Diff}^j(M_2)$ denotes the group of all 
$C^j$-diffeomorphisms of $M_2$, and $S^0(M_i)$ ($i=1,2$) means $M_i$.  
It is easy to show that this definition is independent of the choice of the natural number $m$ with 
$m\geq\max\{m_1,m_2\}$.  

Set $\tau_{\rm min}^k(M,g):=\inf\{\|f\|_{C^k}\,|\,f\in{\rm Emb}_I^k((M,g),\mathbb E^{m(M,g)})\}$, 
where $\|\cdot\|_{C^k}$ denotes the $C^k$-norm of the vector space 
$C^k(M,\mathbb E^m)$ of all $C^k$-maps from $M$ to $\mathbb E^{m(M,g)}$ 
(which is regarded as a vector space) with respect to $g$ and $g_{\mathbb E}$.  
If there exists $f\in{\rm Emb}_I((M,g),\mathbb E^{m(M,g)})$ satisfying $\|f\|_{C^k}=\tau_{\rm min}^k(M,g)$, then 
we call $f$ an {\it ideal $C^k$-isometric embedding} {\it of} $(M,g)$ {\it into} $\mathbb E^{m(M,g)}$.  

\vspace{0.5truecm}

\noindent
{\it Remark 3.1.} If the rank of the second fundamental form of $n(\geq3)$-dimensional compact Riemannian 
hypersurface $(M,g)$ in the Euclidean space $\mathbb E^{n+1}$ is greater than two, then the isometric embedding 
is rigid (see \cite{DR}).  
That is, $\sharp({\rm Emb}_I((M,g),\mathbb E^{n+1})/\equiv)=1$ and hence all elements of 
${\rm Emb}_I((M,g),\mathbb E^{n+1})$ are ideal isometric embeddings, where $\equiv$ denotes the congruence relation 
and $\sharp(\cdot)$ denotes the cardinal number of $(\cdot)$.  

\vspace{0.5truecm}

The following fact holds for $d_{GH}^k|_{(\mathcal{RM}^k_c)_I\times(\mathcal{RM}^k_c)_I}$.  

\vspace{0.5truecm}

\noindent
{\bf Theorem 3.1.} {\sl $d_{GH}^k|_{\mathcal{RM}^k_c\times\mathcal{RM}^k_c}$ is a distance function over 
$\mathcal{RM}^k_c$.} 

\vspace{0.5truecm}

To prove Thoerem 3.1, we prepare the following lemma.  

\vspace{0.5truecm}

\noindent
{\bf Lemma 3.2.} {\sl Let $[(M_i,g_i)]\in\mathcal{RM}_c^k$ ($i=1,2$), $m_1,m_2,m$ be as above and 
$\{f_1^l\}_{l=1}^{\infty}$ be a sequence in ${\rm Emb}_I^{k+1}((M_1,g_1),\mathbb E^{m_1})$ satisfying 
the following conditions:

\vspace{0.15truecm}

(i) The barycenters of $f_1^l(M_1)$ ($l\in\mathbb N$) coincide with one another;

(ii) $\lim\limits_{l\to\infty}\|f_1^l\|_{C^{k+1}}=\tau_{\rm min}^{k+1}(M_1,g_1)$.

\vspace{0.15truecm}

\noindent
Then there exists a sequence $\{f_2^l\}_{l=1}^{\infty}$ in ${\rm Emb}_I^{k+1}((M_2,g_2),\mathbb E^{m_2})$ satisfying 
the following conditions:

\vspace{0.15truecm}

(i') The barycenters of $f_2^l(M_2)$ ($l\in\mathbb N$) coincide with one another;

(ii') $\lim\limits_{l\to\infty}\|f_2^l\|_{C^{k+1}}=\tau_{\rm min}^{k+1}(M_2,g_2)$;

(iii') $\displaystyle{\begin{array}{l}
\displaystyle{\sum_{j=0}^{k+1}\,d_{H,(T^j\mathbb E^m,(g_{\mathbb E})_S^j)}(d^j(\iota_1\circ f_1^l)(S^jM_1),
d^j(\iota_2\circ f_2^l\circ\psi)(S^jM_2))}\\
\displaystyle{\longrightarrow\,\,d_{GH}^k([(M_1,g_1)],[(M_2,g_2)])\,\,\,\,\,(l\to\infty)}
\end{array}}$ 

\noindent
holds for suitable totally geodesic embeddings $\iota_i$'s of $\mathbb E^{m_i}$ into 
$\mathbb E^m$ and a suitable $C^{k+1}$-diffeomorphism $\psi$ of $M_2$.}

\vspace{0.5truecm}

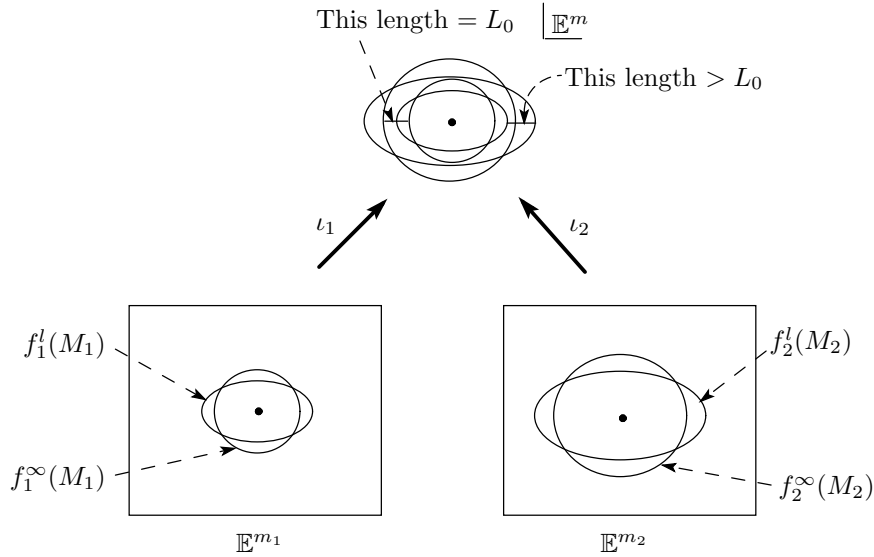
\centerline{
\unitlength 0.1in
\begin{picture}( 60.4000, 34.5000)(-16.2000,-39.9000)
%
\special{pn 8}%
\special{ar 1692 2690 290 160  0.0000000 6.2831853}%
%
\special{pn 8}%
\special{ar 1692 2690 224 218  0.0000000 6.2831853}%
%
\special{pn 8}%
\special{pa 2982 2136}%
\special{pa 2982 3234}%
\special{pa 4302 3234}%
\special{pa 4302 2136}%
\special{pa 4302 2136}%
\special{pa 2982 2136}%
\special{fp}%
%
\special{pn 8}%
\special{ar 3592 2712 448 232  0.0000000 6.2831853}%
%
\special{pn 8}%
\special{ar 3592 2712 348 320  0.0000000 6.2831853}%
%
\special{pn 8}%
\special{pa 1022 2136}%
\special{pa 1022 3234}%
\special{pa 2342 3234}%
\special{pa 2342 2136}%
\special{pa 2342 2136}%
\special{pa 1022 2136}%
\special{fp}%
%
\special{pn 20}%
\special{sh 1}%
\special{ar 3606 2726 10 10 0  6.28318530717959E+0000}%
\special{sh 1}%
\special{ar 3606 2726 10 10 0  6.28318530717959E+0000}%
%
\special{pn 20}%
\special{sh 1}%
\special{ar 1700 2690 10 10 0  6.28318530717959E+0000}%
\special{sh 1}%
\special{ar 1700 2690 10 10 0  6.28318530717959E+0000}%
%
\special{pn 8}%
\special{ar 2712 1168 290 158  0.0000000 6.2831853}%
%
\special{pn 8}%
\special{ar 2712 1168 224 218  0.0000000 6.2831853}%
%
\special{pn 8}%
\special{ar 2700 1170 448 232  0.0000000 6.2831853}%
%
\special{pn 8}%
\special{ar 2698 1164 346 320  0.0000000 6.2831853}%
%
\special{pn 20}%
\special{sh 1}%
\special{ar 2712 1176 10 10 0  6.28318530717959E+0000}%
\special{sh 1}%
\special{ar 2712 1176 10 10 0  6.28318530717959E+0000}%
%
\special{pn 20}%
\special{pa 2016 1934}%
\special{pa 2350 1600}%
\special{fp}%
\special{sh 1}%
\special{pa 2350 1600}%
\special{pa 2290 1632}%
\special{pa 2312 1638}%
\special{pa 2318 1660}%
\special{pa 2350 1600}%
\special{fp}%
%
\special{pn 20}%
\special{pa 3408 1960}%
\special{pa 3082 1592}%
\special{fp}%
\special{sh 1}%
\special{pa 3082 1592}%
\special{pa 3110 1654}%
\special{pa 3116 1632}%
\special{pa 3140 1628}%
\special{pa 3082 1592}%
\special{fp}%
%
\special{pn 8}%
\special{pa 3200 738}%
\special{pa 3390 738}%
\special{fp}%
%
\special{pn 8}%
\special{pa 3190 738}%
\special{pa 3190 552}%
\special{fp}%
\put(32.3500,-7.1000){\makebox(0,0)[lb]{$\mathbb E^m$}}%
\put(15.8200,-33.3100){\makebox(0,0)[lt]{$\mathbb E^{m_1}$}}%
\put(34.8800,-33.3100){\makebox(0,0)[lt]{$\mathbb E^{m_2}$}}%
\put(21.1000,-17.6000){\makebox(0,0)[rb]{$\iota_1$}}%
\put(33.2500,-17.6600){\makebox(0,0)[lb]{$\iota_2$}}%
\put(9.0000,-24.1000){\makebox(0,0)[rb]{{\small $f_1^l(M_1)$}}}%
\put(9.0000,-29.7000){\makebox(0,0)[rt]{{\small $f_1^{\infty}(M_1)$}}}%
\put(43.7300,-23.9900){\makebox(0,0)[lb]{{\small $f_2^l(M_2)$}}}%
\put(44.2000,-30.2000){\makebox(0,0)[lt]{{\small $f_2^{\infty}(M_2)$}}}%
%
\special{pn 8}%
\special{pa 950 3024}%
\special{pa 1574 2884}%
\special{da 0.070}%
\special{sh 1}%
\special{pa 1574 2884}%
\special{pa 1504 2878}%
\special{pa 1522 2896}%
\special{pa 1512 2918}%
\special{pa 1574 2884}%
\special{fp}%
%
\special{pn 8}%
\special{pa 4382 3076}%
\special{pa 3814 2972}%
\special{da 0.070}%
\special{sh 1}%
\special{pa 3814 2972}%
\special{pa 3876 3004}%
\special{pa 3866 2982}%
\special{pa 3882 2964}%
\special{pa 3814 2972}%
\special{fp}%
%
\special{pn 8}%
\special{pa 4356 2374}%
\special{pa 4012 2638}%
\special{da 0.070}%
\special{sh 1}%
\special{pa 4012 2638}%
\special{pa 4078 2612}%
\special{pa 4054 2604}%
\special{pa 4054 2580}%
\special{pa 4012 2638}%
\special{fp}%
%
\special{pn 8}%
\special{pa 950 2364}%
\special{pa 1420 2628}%
\special{da 0.070}%
\special{sh 1}%
\special{pa 1420 2628}%
\special{pa 1372 2578}%
\special{pa 1374 2602}%
\special{pa 1352 2614}%
\special{pa 1420 2628}%
\special{fp}%
\put(25.0000,-39.9000){\makebox(0,0)[lt]{$\,$}}%
%
\special{pn 8}%
\special{pa 2360 1170}%
\special{pa 2480 1170}%
\special{fp}%
\put(10.5000,-36.2000){\makebox(0,0)[lt]{{\small $f_i^{\infty}:=\lim\limits_{l\to\infty}f_i^l$ ($i=1,2$) (Assume that these limits exist.)}}}%
%
\special{pn 8}%
\special{pa 2230 750}%
\special{pa 2400 1170}%
\special{da 0.070}%
\special{sh 1}%
\special{pa 2400 1170}%
\special{pa 2394 1102}%
\special{pa 2380 1122}%
\special{pa 2356 1116}%
\special{pa 2400 1170}%
\special{fp}%
\put(20.0000,-7.1000){\makebox(0,0)[lb]{{\small This length $=L_0$}}}%
\put(33.0000,-10.1000){\makebox(0,0)[lb]{{\small This length $>L_0$}}}%
%
\special{pn 8}%
\special{ar 3290 1120 200 190  3.1036798 3.4113721}%
\special{ar 3290 1120 200 190  3.5959875 3.9036798}%
\special{ar 3290 1120 200 190  4.0882952 4.3959875}%
\special{ar 3290 1120 200 190  4.5806029 4.6827241}%
%
\special{pn 8}%
\special{pa 3090 1120}%
\special{pa 3090 1180}%
\special{da 0.070}%
\special{sh 1}%
\special{pa 3090 1180}%
\special{pa 3110 1114}%
\special{pa 3090 1128}%
\special{pa 3070 1114}%
\special{pa 3090 1180}%
\special{fp}%
%
\special{pn 8}%
\special{pa 3000 1180}%
\special{pa 3150 1180}%
\special{fp}%
\put(10.5000,-38.3000){\makebox(0,0)[lt]{{\small ($f_i^{\infty}\,:\,$an ideal $C^{k+1}$-isometric embedding)}}}%
\end{picture}%
\hspace{6.5truecm}}

\vspace{0.2truecm}

\centerline{{\bf Figure 6$\,\,\,:\,\,$ About the statement of Lemma 3.2}}

\vspace{0.5truecm}

It is clear that the distance $d_{GH}^k([(M_1,g_1)],[(M_2,g_2)])$ is attained when the $C^{k+1}$-norms of 
isometric embeddings of $(M_i,g_i)$'s ($i=1,2$) into $\mathbb E^{m_i}$ are as small as possible (i.e., 
the isometric embeddings are as close as to ideal isometric embeddings possible), 
$\mathbb E^{m_i}$ ($i=1,2$) are isometrically embedded into $\mathbb E^m$ by suitable totally geodesic 
isometric embeddings, where we need to replace $g_2$ to another Riemannian metric belonging to the isometric class 
$[g_2]$ of $g_2$.  Hence it is clear that the statement of this lemma holds.  

As an example, we shall consider the case of $(M_i,g_i)=\mathbb S^n(r_i^{-2})\,\,\,(i=1,2)$.  
In this case, we can confirm that the statement of Lemma 3.2 holds as follows.  
Since $\mathbb S^n(r_i^{-2})$'s are isometrically embedded into 
$\mathbb E^{n+1}$, $m_i$ ($i=1,2$) in Lemma 3.2 are equal to $n+1$ and hence $m$ in Lemma 3.2 also is equal to 
$n+1$.  First we consider the case of $n\geq 2$.  In this case, an $C^{k+1}$-isometric embedding of 
$\mathbb S^n(r_i^{-2})$ into $\mathbb E^{n+1}$ is unique up to congruence and it is a totally umbilic 
$C^{\infty}$-isometric embedding.  Naturally it is an ideal isometric embedding.  
If two totally umbilic isometric embeddings of $\mathbb S^n(r_i^{-2})$ into 
$\mathbb E^{n+1}$ have the same baryccenter, their images coincide with each other.  
Let $f^u_i$ ($i=1,2$) be a totally umbilic isometric embeddding of $\mathbb S^n(r_i^{-2})$ into 
$\mathbb E^{n+1}$ whose barycenter is equal to the origin $o$ of $\mathbb E^{n+1}$.  
Let $\{f_1^l\}_{l=1}^{\infty}$ be a sequence in ${\rm Emb}_I(\mathbb S^n(r_i^{-2}),\mathbb E^{n+1})$ having 
the same barycenter.  Denote by $p_G$ the barycenter.  
Then each $f_1^l$ is expressed as $f_1^l=\phi_l\circ f^1_1$ for some isometry $\phi_l$ of $\mathbb E^{n+1}$ fixing 
$p_G$.  Define $f_2^l\in{\rm Emb}_I(\mathbb S^n(r_2^{-2}),\mathbb E^{n+1})$ 
$$f_2^l:=\phi_l\circ\frac{r_2}{r_1}f_1^1\circ\frac{r_1}{r_2}\,{\rm id}_{\mathbb E^{n+1}}|_{\mathbb S^n(r_2^{-2})}.$$
Let $\iota_1=\iota_2={\rm id}_{\mathbb E^{n+1}}$.  Then we have 
$$\begin{array}{r}
\displaystyle{\sum_{j=0}^{k+1}d_{H,(T^j\mathbb E^{n+1},(g_{\mathbb E})_S^j)}
(d^j(\iota_1\circ f_1^l)(S^j(\mathbb S^n(r_1^{-2})),d^j(\iota_2\circ f_2^l)(S^j(\mathbb S^n(r_2^{-2})))}\\
\displaystyle{=d_{GH}^k(\mathbb S^n(r_1^{-2}),\mathbb S^n((r_2^{-2})).}
\end{array}$$

Next we consider the case of $n=1$.  Let $f_1^u$ be an isometric embedding of $S^1(r_1^{-2})$ into 
$\mathbb E^2$ whose image is equal to the circle of radius $r_1$ centered at $o$.  
It is clear that $f_1^u$ is an ideal $C^{k+1}$-isometric embedding of $\mathbb S^1(r_1^{-2})$ into $\mathbb E^2$.  
Let $\{f_1^l\}_{l=1}^{\infty}$ be a sequence in ${\rm Emb}_I(\mathbb S^1(r_1^{-2}),\mathbb E^2)$ 
having the same barycenter and satisfying 
$\lim\limits_{l\to\infty}\|f_1^l\|_{C^{k+1}}=\tau_{\rm min}^{k+1}(\mathbb S^1(r_1^{-2}))$.  
Let $p_G$ be the same barycenter of $f_1^l(S^1(r_1^{-2}))$'s and $\tau_{\pm p_G}:\mathbb E^2\to\mathbb E^2$ be 
the parallel translation by $\pm p_G$ (i.e., $\tau_{\pm p_G}(p):=p\pm p_G\,\,\,(p\in\mathbb E^2)$).  
Then, it follows from 
$\|f_1^l\|_{C^{k+1}}\to\tau_{\rm min}^{k+1}(\mathbb S^1(r_1^{-2}))$ that 
$\|f_1^l-\phi_l\circ\tau_{p_G}\circ f_1^u\|_{C^{k+1}}\to 0\,\,\,(l\to\infty)$ holds for a family 
$\{\phi_l\}_{l=1}^{\infty}$ of isometries of $\mathbb E^2$ fixing $p_G$.  
Set 
$$f_2^l:=\phi_l\circ\tau_{p_G}\circ\frac{r_2}{r_1}\cdot{\rm id}_{\mathbb E^2}\circ\tau_{-p_G}\circ f_1^l\circ
\frac{r_1}{r_2}\cdot{\rm id}_{\mathbb E^2}|_{S^1(r_2^{-2})}$$
and let $\iota_1=\iota_2={\rm id}_{\mathbb E^2}$.  Then we have 
$$\begin{array}{r}
\displaystyle{\sum_{j=0}^{k+1}d_{H,(T^j\mathbb E^2,(g_{\mathbb E})_S^j)}(d^j(\iota_1\circ f_1^l)(S^j(\mathbb S^1(r_1^{-2})),
d^j(\iota_2\circ f_2^l)(S^j(\mathbb S^1(r_2^{-2})))}\\
\to d_{GH}^k([\mathbb S^1(r_1^{-2})],[\mathbb S^1(r_2^{-2})])\quad\,\,(l\to\infty).
\end{array}$$
(see Figure 7).  

By using this lemma, we prove Theorem 3.1.  

\vspace{0.5truecm}

\noindent
{\it Proof of Theorem 3.1.} 
holds for suitable totally geodesic embeddings $\iota_i$'s ($1=1,2$) of $\mathbb E^{m_i}$ into 
$\mathbb E^{m_{12}}$.  
First we show that $d^k_{GH}$ is a pseudo-distance function.  
It is clear that $d^k_{GH}$ satisfies the conditions other than the triangle inequality 
in the definition of the pseudo-distance function.  We show that it satisfies the triangle inequality.  
Take $[(M_i,g_i)]\in\mathcal{RM}_c^k$ ($i=1,2,3$).  

\newpage


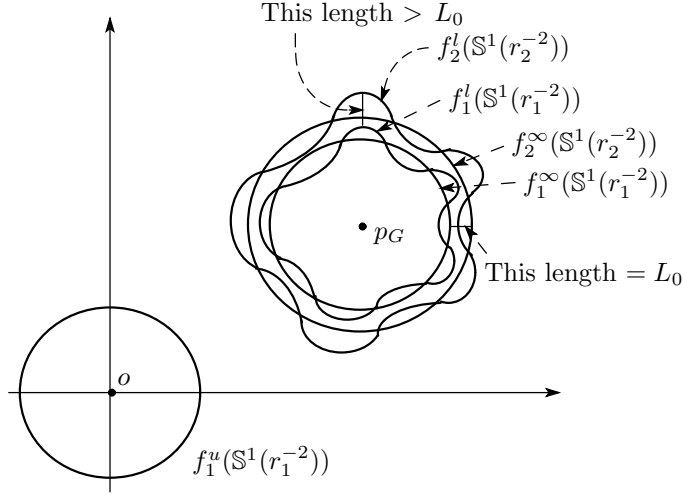
\centerline{
\unitlength 0.1in
\begin{picture}( 28.8000, 31.3000)( 10.8000,-35.5000)
\put(15.4000,-35.5000){\makebox(0,0)[lt]{$\,$}}%
%
\special{pn 13}%
\special{ar 2918 1630 586 560  0.0000000 6.2831853}%
%
\special{pn 20}%
\special{sh 1}%
\special{ar 1620 2510 10 10 0  6.28318530717959E+0000}%
\special{sh 1}%
\special{ar 1620 2510 10 10 0  6.28318530717959E+0000}%
%
\special{pn 13}%
\special{ar 2932 1160 170 220  3.5102886 5.8175397}%
%
\special{pn 13}%
\special{ar 3312 940 250 310  1.2570314 2.7589330}%
%
\special{pn 13}%
\special{ar 3372 1370 190 130  4.7123890 6.2831853}%
\special{ar 3372 1370 190 130  0.0000000 0.7940186}%
%
\special{pn 13}%
\special{ar 3612 1630 180 200  2.2142974 4.1139701}%
%
\special{pn 13}%
\special{ar 3382 1880 140 160  5.6580178 6.2831853}%
\special{ar 3382 1880 140 160  0.0000000 1.8640993}%
%
\special{pn 13}%
\special{ar 3362 2200 350 160  3.0658756 4.6509550}%
%
\special{pn 13}%
\special{ar 2822 2140 210 160  0.4589992 3.1848989}%
%
\special{pn 13}%
\special{ar 2362 2240 260 380  4.7570686 6.0293236}%
%
\special{pn 13}%
\special{ar 2442 1610 200 260  1.8908885 4.5714850}%
%
\special{pn 13}%
\special{ar 2202 990 590 390  0.1680014 1.2070037}%
\put(16.5000,-24.7000){\makebox(0,0)[lb]{$o$}}%
%
\special{pn 8}%
\special{pa 1080 2510}%
\special{pa 3960 2510}%
\special{fp}%
\special{sh 1}%
\special{pa 3960 2510}%
\special{pa 3894 2490}%
\special{pa 3908 2510}%
\special{pa 3894 2530}%
\special{pa 3960 2510}%
\special{fp}%
%
\special{pn 8}%
\special{pa 1610 3030}%
\special{pa 1610 550}%
\special{fp}%
\special{sh 1}%
\special{pa 1610 550}%
\special{pa 1590 618}%
\special{pa 1610 604}%
\special{pa 1630 618}%
\special{pa 1610 550}%
\special{fp}%
%
\special{pn 13}%
\special{ar 2918 1630 472 446  0.0000000 6.2831853}%
%
\special{pn 13}%
\special{ar 2936 1284 134 164  3.5110885 5.8157354}%
%
\special{pn 13}%
\special{ar 3236 1120 198 230  1.2569020 2.7583525}%
%
\special{pn 13}%
\special{ar 3284 1440 150 96  4.7123890 6.2831853}%
\special{ar 3284 1440 150 96  0.0000000 0.7936624}%
%
\special{pn 13}%
\special{ar 3472 1632 142 148  2.2103560 4.1170498}%
%
\special{pn 13}%
\special{ar 3302 1820 110 120  5.6603276 6.2831853}%
\special{ar 3302 1820 110 120  0.0000000 1.8653867}%
%
\special{pn 13}%
\special{ar 3282 2060 276 120  3.0608061 4.6508765}%
%
\special{pn 13}%
\special{ar 2852 2010 166 118  0.4607021 3.1860549}%
%
\special{pn 13}%
\special{ar 2482 2090 206 282  4.7576411 6.0281426}%
%
\special{pn 13}%
\special{ar 2552 1620 158 194  1.8893857 4.5721740}%
%
\special{pn 13}%
\special{ar 2362 1158 464 290  0.1682814 1.2064016}%
\put(33.8000,-10.5000){\makebox(0,0)[lb]{$f_1^l(\mathbb S^1(r_1^{-2}))$}}%
%
\special{pn 20}%
\special{sh 1}%
\special{ar 2932 1640 10 10 0  6.28318530717959E+0000}%
\special{sh 1}%
\special{ar 2932 1640 10 10 0  6.28318530717959E+0000}%
\put(29.9200,-17.3000){\makebox(0,0)[lb]{$p_G$}}%
%
\special{pn 13}%
\special{ar 1610 2510 472 446  0.0000000 6.2831853}%
\put(20.4000,-27.9000){\makebox(0,0)[lt]{$f_1^u(\mathbb S^1(r_1^{-2}))$}}%
%
\special{pn 8}%
\special{pa 3336 980}%
\special{pa 3010 1140}%
\special{da 0.070}%
\special{sh 1}%
\special{pa 3010 1140}%
\special{pa 3080 1128}%
\special{pa 3058 1116}%
\special{pa 3062 1094}%
\special{pa 3010 1140}%
\special{fp}%
\put(33.2000,-7.8000){\makebox(0,0)[lb]{$f_2^l(\mathbb S^1(r_2^{-2}))$}}%
\put(37.7000,-14.9000){\makebox(0,0)[lb]{$f_1^{\infty}(\mathbb S^1(r_1^{-2}))$}}%
\put(37.1000,-12.8000){\makebox(0,0)[lb]{$f_2^{\infty}(\mathbb S^1(r_2^{-2}))$}}%
%
\special{pn 8}%
\special{pa 2932 1110}%
\special{pa 2932 940}%
\special{fp}%
%
\special{pn 8}%
\special{ar 2850 730 268 300  1.5707963 1.7820639}%
\special{ar 2850 730 268 300  1.9088245 2.1200921}%
\special{ar 2850 730 268 300  2.2468527 2.4581203}%
\special{ar 2850 730 268 300  2.5848808 2.7961484}%
\special{ar 2850 730 268 300  2.9229090 3.1341766}%
%
\special{pn 8}%
\special{pa 2862 1030}%
\special{pa 2932 1030}%
\special{da 0.070}%
\special{sh 1}%
\special{pa 2932 1030}%
\special{pa 2866 1010}%
\special{pa 2880 1030}%
\special{pa 2866 1050}%
\special{pa 2932 1030}%
\special{fp}%
\put(24.0000,-5.9000){\makebox(0,0)[lb]{{\small This length $>\,L_0$}}}%
\put(35.7000,-18.3000){\makebox(0,0)[lt]{{\small This length $=L_0$}}}%
\put(15.0000,-31.4000){\makebox(0,0)[lt]{$f_i^{\infty}:=\lim\limits_{l\to\infty}f_i^l$ ($i=1,2$) {\small (Assume that these limits exist.)}}}%
%
\special{pn 8}%
\special{ar 3280 980 240 276  3.3422002 3.5752099}%
\special{ar 3280 980 240 276  3.7150157 3.9480254}%
\special{ar 3280 980 240 276  4.0878313 4.3208410}%
\special{ar 3280 980 240 276  4.4606468 4.6936565}%
%
\special{pn 8}%
\special{pa 3040 950}%
\special{pa 3030 980}%
\special{da 0.070}%
\special{sh 1}%
\special{pa 3030 980}%
\special{pa 3070 924}%
\special{pa 3048 930}%
\special{pa 3032 910}%
\special{pa 3030 980}%
\special{fp}%
%
\special{pn 8}%
\special{pa 3670 1220}%
\special{pa 3410 1320}%
\special{da 0.070}%
\special{sh 1}%
\special{pa 3410 1320}%
\special{pa 3480 1316}%
\special{pa 3460 1302}%
\special{pa 3466 1278}%
\special{pa 3410 1320}%
\special{fp}%
%
\special{pn 8}%
\special{pa 3730 1410}%
\special{pa 3340 1440}%
\special{da 0.070}%
\special{sh 1}%
\special{pa 3340 1440}%
\special{pa 3408 1456}%
\special{pa 3394 1436}%
\special{pa 3406 1416}%
\special{pa 3340 1440}%
\special{fp}%
%
\special{pn 8}%
\special{pa 2580 630}%
\special{pa 2580 730}%
\special{da 0.070}%
%
\special{pn 8}%
\special{pa 3390 1640}%
\special{pa 3500 1640}%
\special{fp}%
%
\special{pn 8}%
\special{pa 3580 1780}%
\special{pa 3470 1640}%
\special{da 0.070}%
\special{sh 1}%
\special{pa 3470 1640}%
\special{pa 3496 1706}%
\special{pa 3504 1682}%
\special{pa 3528 1680}%
\special{pa 3470 1640}%
\special{fp}%
\put(15.0000,-33.6000){\makebox(0,0)[lt]{{\small ($f_i^{\infty}\,:\,$ ideal $C^{\infty}$-isometric embedding)}}}%
\end{picture}%
\hspace{3truecm}}

\vspace{0.2truecm}

\centerline{{\bf Figure 7$\,\,\,:\,\,$ An example of sequences as in Lemma 3.2}}

\vspace{0.5truecm}

\noindent
Set $a_{ij}:=d_{GH}^k([(M_i,g_i)],[(M_j,g_j)])$.  
Set $m_{123}:=\max\{m_1,m_2,m_3\}$.  
Take a sequence $\{f_2^l\}_{l=1}^{\infty}$ in ${\rm Emb}_I((M_2,g_2),\mathbb E^{m_2})$ having the same barycenter 
and satisfying $\|f_2^l\|_{C^{k+1}}\to\tau_{\rm min}^{k+1}(M_2,g_2)\,\,\,\,(l\to\infty)$.  
Then, according to Lemma 3.2, there exist sequences $\{f_i^l\}_{l=1}^{\infty}$ ($i=1,3$) in 
${\rm Emb}_I^{k+1}((M_i,g_i),\mathbb E^{m_i})$ having the same barycenter and satisfying the following conditions:

\vspace{0.15truecm}

(i)\ $\|f_i^l\|_{C^{k+1}}\to\tau_{\rm min}^{k+1}(M_i,g_i)\,\,\,\,(l\to\infty);$ 

(ii) $\displaystyle{\begin{array}{r}
\displaystyle{\sum_{j=0}^{k+1}\,d_{H,(T^j\mathbb E^m,(g_{\mathbb E})_S^j)}
(d^j(\iota_{2i}\circ f_2^l)(S^jM_2),d^j(\iota_i\circ f_i^l\circ\psi_i)(S^jM_i))\to a_{2i}}\\
\displaystyle{(l\to\infty)}
\end{array}}$\newline
($i=1,3$) holds for suitable totally geodesic embeddings $\iota_{2i}$ ($i=1,3$) of $\mathbb E^{m_2}$ into 
$\mathbb E^{m_{123}}$, suitable totally geodesic embeddings $\iota_i$'s ($i=1,3$) of $\mathbb E^{m_i}$ into 
$\mathbb E^{m_{123}}$ and a suitable $C^{k+1}$-diffeomorphism $\psi_i$ of $M_i$.  

\vspace{0.15truecm}

\noindent
Then we have 
\begin{align*}
&d_{H,(T^j\mathbb E^{m_{123}},(g_{\mathbb E})_S^j)}(d^j(\iota_{21}^{-1}\circ\iota_1\circ f_1^l\circ\psi_1)(S^j(M_1)),
d^j(\iota_{23}^{-1}\circ\iota_3\circ f_3^l\circ\psi_3)(S^j(M_3)))\\
\leq&d_{H,(T^j\mathbb E^{m_{123}},(g_{\mathbb E})_S^j)}
(d^j(\iota_{21}^{-1}\circ\iota_1\circ f_1^l\circ\psi_1)(S^j(M_1)),d^jf_2^l(S^j(M_2)))\\
&+d_{H,(T^j\mathbb E^{m_{123}},(g_{\mathbb E})_S^j)}(d^jf_2^l(S^j(M_2)),
d^j(\iota_{23}^{-1}\circ\iota_3\circ f_3^l\circ\psi_3)(S^j(M_3)))
\end{align*}
and hence 
{\small
$$\lim_{l\to\infty}\,\sum_{j=0}^{k+1}\,
d_{H,(T^j\mathbb E^{m_{123}},(g_{\mathbb E})_S^j)}(d^j(\iota_{21}^{-1}\circ\iota_1\circ f_1^l\circ\psi_1)(S^j(M_1)),
d^j(\iota_{23}^{-1}\circ\iota_3\circ f_3^l\circ\psi_3)(S^j(M_3)))$$
}
is smaller than or equal to $a_{12}+a_{23}$.
Therefore we obtain $a_{13}\leq a_{12}+a_{23}$.  Thus $d^k_{GH}$ satisfies the triangle inequality and hence 
it is a pseudo-distance function.  

Furthermore we show that $d^k_{GH}$ is a distance function.  
Assume that $d^k_{GH}([(M_1,g_1)],$\newline
$[(M_2,g_2)])=0$.  
Set $m_{12}:=\max\{m_1,m_2\}$.  
Take a sequence $\{f_1^l\}_{l=1}^{\infty}$ in ${\rm Emb}_I^{k+1}((M_1,g_1),\mathbb E^{m_1})$ having the same 
barycenter and satisfying $\|f_1^l\|_{C^{k+1}}\to\tau_{\min}^{k+1}(M_1,g_1)\,\,\,(l\to\infty)$.  
Then, according to Lemma 3.2, there exist a sequence $\{f_2^l\}_{l=1}^{\infty}$ in 
${\rm Emb}_I^{k+1}((M_2,g_2),\mathbb E^{m_2})$ having the same barycenter and satisfying the following conditions:

\vspace{0.15truecm}

(i)\ $\|f_2^l\|_{C^{k+1}}\to\tau_{\rm min}^{k+1}(M_i,g_i)\,\,\,\,(l\to\infty);$ 

(ii) $\displaystyle{\begin{array}{l}
\displaystyle{\lim_{l\to\infty}\,\sum_{j=0}^{k+1}\,
d_{H,(T^j\mathbb E^{m_{12}},(g_{\mathbb E})_S^j)}
(d^j(\iota_1\circ f_1^l)(S^jM_1),d^j(\iota_2\circ f_2^l\circ\psi_2)(S^jM_2))}\\
\displaystyle{=d^k_{GH}([(M_1,g_1)],[(M_2,g_2)])(=0)}
\end{array}
}$

\vspace{0.15truecm}

\noindent
holds for suitable totally geodesic embeddings $\iota_i$'s ($1=1,2$) of $\mathbb E^{m_i}$ into 
$\mathbb E^{m_{12}}$ and a suitable $C^{k+1}$-diffeomorphism $\psi_2$ of $M_2$.  

\vspace{0.15truecm}

\noindent
This implies that $M_1$ and $M_2$ are $C^{k+1}$-diffeomorphic.  
Take a $C^{k+1}$-diffeomorphism $\psi$ of $M_1$ onto $M_2$.  Since $[(M_1,\psi^{\ast}g_2)]=[(M_2,g_2)]$, 
we have $d^k_{GH}([(M_1,g_1)],[(M_1,\psi^{\ast}g_2)])=0$ by the assumption.  
According to Lemma 3.2, there exists a sequence $\{\widehat f_1^l\}_{l=1}^{\infty}$ in 
${\rm Emb}_I^{k+1}((M_1,\psi^{\ast}g_2),$\newline
$\mathbb E^{m_2})$ having the same barycenter and satisfying the following conditions:

\vspace{0.15truecm}

(i)\ $\|\widehat f_1^l\|_{C^{k+1}}\to\tau_{\rm min}^{k+1}(M_1,\psi^{\ast}g_2)\,\,\,\,(l\to\infty);$ 

(ii) $\displaystyle{\begin{array}{l}
\displaystyle{\lim_{l\to\infty}\,\sum_{j=0}^{k+1}\,d_{H,(T^j\mathbb E^{m_{12}},(g_{\mathbb E})_S^j)}
(d^j(\iota_1\circ f_1^l)(S^jM_1),d^j(\widehat{\iota}_1\circ\widehat f_1^l\circ\widehat{\psi}_1)
(\widehat S^jM_1))}\\
\displaystyle{=d^k_{GH}([(M_1,g_1)],[(M_1,\psi^{\ast}g_2)])(=0)}
\end{array}}$\newline
holds for a suitable totally geodesic isometric embedding $\widehat{\iota}_1$ of $(M_1,\psi^{\ast}g_2)$ 
into $\mathbb E^{m_2}$ and a suitable $C^{k+1}$-diffeomorphism $\widehat{\psi}_1$ of $M_1$, where $S^{k+1}M_1$ and 
$\widehat S^{k+1}M_1$ denote the unit tangent bundles of 
the Riemannian manifolds $(T^{k+1}M_1,(g_1)_S^k)$ and $(T^{k+1}M_1,(\psi^{\ast}g_2)_S^k)$, respectively.  

\vspace{0.15truecm}

\noindent
From the condition (ii), we can show that 
$\lim\limits_{l\to\infty}\|(\iota_1\circ f_1^l)
-(\widehat{\iota}_1\circ\widehat f_1^l\circ\widehat{\psi}_1)\|_{C^{k+1}}=0$ holds.  Then we have 
$\lim\limits_{l\to\infty}\|(\iota_1\circ f_1^l)^{\ast}g_{\mathbb E}
-(\widehat{\iota}_1\circ\widehat f_1^l\circ\widehat{\psi}_1)^{\ast}g_{\mathbb E}\|_{C^k}=0$, 
where $\|\cdot\|_{C^k}$ denotes the $C^k$-norm of the space of all $C^k$-sections of the tensor bundle 
$T^{\ast}M_1\otimes T^{\ast}M_1$.  
Therefore, by noticing $(\iota_1\circ f_1^j)^{\ast}g_{\mathbb E}=g_1$ and 
$(\widehat{\iota}_1\circ\widehat f_1^j\circ\widehat{\psi}_1)^{\ast}g_{\mathbb E}
=(\psi\circ\widehat{\psi}_1)^{\ast}g_2$, 
we obtain $g_1=(\psi\circ\widehat{\psi}_1)^{\ast}g_2$ and hence $[(M_1,g_1)]=[(M_1,g_2)]$.  
Therefore $d_{GH}^k$ is a distance function over $\mathcal{RM}^k_c$.  \qed

\vspace{0.5truecm}

\noindent
{\bf Problem.}\ {\sl Does $d^k_{GH}$ coincide with $d_{GH}|_{\mathcal{RM}_c^k\times\mathcal{RM}_c^k}$?}

\vspace{0.5truecm}

If this problem were solved affirmatively, then we can define a completion of $(\mathcal{RM}_c^k,d_{GH}^k)$ as 
$(\overline{\mathcal{RM}_c^k},d_{GH}|_{\overline{\mathcal{RM}_c^k}\times\overline{\mathcal{RM}_c^k}})$, where 
$\overline{\mathcal{RM}_c^k}$ denotes the closure of $\mathcal{RM}_c^k$ in $(\mathcal M_c,d_{GH})$.  

\section{Hamilton convergence} 
R. S. Hamilton (\cite{H}) introduced the convergence of a sequence of complete marked Riemannian manifolds.  
We recall the definition of this convergence.  Let $(M,g)$ be a complete Riemannian manifold, $p$ be a point of 
$M$ and $O$ be an orthonormal frame of $(M,g)$ at $p$.  Then $(M,g,p,O)$ is called a {\it complete marked 
Riemannian manifold}.  Let $(M,g,p,O)$ be a complete marked Riemannian manifold and 
$\{(M_i,g_i,p_i,O_i)\}_{i=1}^{\infty}$ be a sequence of complete marked Riemannan manifolds.  
Assume that there exist a sequence $\{U_i\}_{i=1}^{\infty}$ of open subsets of $M$ and 
a sequence $\{V_i\}_{i=1}^{\infty}$ of open subsets of $M_i$ and a sequence $\{\phi_i:U_i\to V_i\}_{i=1}^{\infty}$ 
of $C^{\infty}$-diffeomorphisms with $\phi_i(p)=p_i$ and $(\phi_i)_{\ast}(O)=O_i$ satisfying 
the following condition:

\vspace{0.15truecm}

\noindent
(H)\ for any compact subset $K$ of $M$, there exists $i_K\in\mathbb N$ such that $K\subset U_i\,\,\,(i\geq i_K)$ 

and that $\{\phi_i^{\ast}g_i|_K\}_{i=i_K}^{\infty}$ converges to $g|_K$ in $C^{\infty}$-topology.  

\vspace{0.15truecm}

\noindent
Then R. S. Hamilton (\cite{H}) called that the sequence $\{(M_i,g_i,p_i,O_i)\}_{i=1}^{\infty}$ converges to 
$(M,g,p,$\newline
$O)$.  In this paper, we call this convergence {\it Hamilton convergence}.  
In particuar, in the case where $M$ and $M_i$ are compact, this convergence may be defined as follows.  
The sequence $\{(M_i,g_i,p_i,O_i)\}_{i=1}^{\infty}$ converges to $(M,g,p,O)$ if and only if 
there exists a sequence $\{\phi_i:M\to M_i\}_{i=1}^{\infty}$ of $C^{\infty}$-diffeomorphisms with 
$\phi_i(p)=p_i$ and $(\phi_i)_{\ast}(O)=O_i$ such that $\{\phi_i^{\ast}g_i\}_{i=1}^{\infty}$ converges to $g$ 
in $C^{\infty}$-topology.  
Furthermore, in this case, we do not need the base points and the base orthonormal frames.  
Hence the Hamilton convergence of a sequence of compact Riemannian manifolds may be defined as follows.  
If there exists a sequence $\{\phi_i:M\to M_i\}_{i=1}^{\infty}$ of $C^{\infty}$-diffeomorphisms such that 
$\{\phi_i^{\ast}g_i\}_{i=1}^{\infty}$ converges to $g$ with respect to $C^{\infty}$-topology, then 
$\{(M_i,g_i)\}_{i=1}^{\infty}$ {\it converges to} $(M,g)$ {\it in the sense of Hamilton}.  
Also, we call then $\{(M_i,g_i)\}_{i=1}^{\infty}$ and the sequence $\{[(M_i,g_i)]\}_{i=1}^{\infty}$ in 
$\mathcal{RM}_c^k$ a {\it Hamilton convergent sequence}.  Also, we define $C^k$-Hamilton covergence as follows.  
If there exists a sequence $\{\phi_i:M\to M_i\}_{i=1}^{\infty}$ of $C^k$-diffeomorphisms such that 
$\{\phi_i^{\ast}g_i\}_{i=1}^{\infty}$ converges to $g$ with respect to $C^k$-norm, then 
$\{(M_i,g_i)\}_{i=1}^{\infty}$ {\it converges to} $(M,g)$ {\it in the sense of Hamilton}.  
Also, we call then the sequence $\{(M_i,g_i)\}_{i=1}^{\infty}$ and the sequence 
$\{[(M_i,g_i)]\}_{i=1}^{\infty}$ in $\mathcal{RM}_c^k$ ``{\it $C^k$-Hamilton convergent sequence}''.  

For the convegence in $d_{GH}^k$ and $C^k$-Hamilton convergence, we prove the following fact.  

\vspace{0.5truecm}

\noindent
{\bf Theorem 4.1.} {\sl For a sequence in $\mathcal{RM}_c^k$, 
the convergence in $d_{GH}^k$ coincides with $C^k$-Hamilton convergence.}

\vspace{0.5truecm}

\noindent
{\it Proof.} Let $\nabla$ be the Riemannian connection of $g$, $\widehat g$ be the fibre metric of the 
$(0,l)$-tensor bundle $T^{(0,l)}M$ of $M$ associated to $g$ and $\widehat{\nabla}$ the connection of 
$T^{(0,l)}M$ associated to $\nabla$.  
The $C^k$-norm $\|\sigma\|_{C^k}$ of a section $\sigma$ of $T^{(0,2)}M$ is defined by 
$$\|\sigma\|_{C^k}:=\sum_{i=0}^k\,\mathop{\sup}_{p\in M}\|(\widehat{\nabla}^i\sigma)_p\|_{\widehat g_p}.
\leqno{(4.1)}$$
Assume that $\{[(M_i,g_i)]\}_{i=1}^{\infty}$ converges to $[(M,g)]$ in $d_{GH}^k$.  
Set $\varepsilon_i:=d_{GH}^k([(M_i,g_i)],[(M,g)])$ and $\widehat m_i:=\max\{m_i,m\}$.  
Then, for each $i$, we can take a sequence $\{f_i^l\}_{l=1}^{\infty}$ in ${\rm Emb}_I^{k+1}((M_i,g_i),$\newline
$E^{m_i})$, a sequence $\{\bar f_i^l\}_{l=1}^{\infty}$ in ${\rm Emb}_I^{k+1}((M,g),E^m)$ and 
$\bar{\psi}_i\in{\rm Diff}^{k+1}(M)$ satisfying 
$$\begin{array}{r}
\displaystyle{\sum_{j=0}^{k+1}\,d_{H,(T^j\mathbb E^{\widehat m_i},d_{(g_{\mathbb E})_S^j})}
\left(d^j(\iota_i\circ f_i^l)(S^jM_i),\,d^j(\bar{\iota}_i\circ\bar f_i^l\circ\bar{\psi}_i)(S^jM)\right)\to
\varepsilon_i}\\
(l\to\infty).
\end{array}
\leqno{(4.2)}$$
Denote by $\iota_i^l$ (resp. $\bar{\iota}_i^l$) be the inclusion map of $f_i^l(M_i)$ (resp. $\bar f_i^l(M)$) 
into $\mathbb E^{\widehat m_i}$.  For a sufficiently large $i_0$, we can show that 
$M_i$ ($i\geq i_0$) are $C^{k+1}$-diffeomorphic to $M$.  
For each $i\geq i_0$, we can take a sequence $\{\psi_i^l\}_{l=1}^{\infty}$ of $C^{k+1}$-diffeomorphisms 
$\psi_i^l$'s of $\bar f_i^l(M)$ onto $f_i^l(M_i)$ such that 
$$\|(\iota_i^l\circ\psi_i^l)^{\ast}g_{\mathbb E}-(\bar{\iota}_i^l\circ\bar{\psi}_i)^{\ast}g_{\mathbb E}\|_{C^k}
\to c_1\varepsilon_i\quad\,\,(l\to\infty),$$
where $c_1$ is a positive constant depending only on $g$.  
Set $\phi_i^l:=(f_i^l)^{-1}\circ\psi_i^l\circ\bar f_i^l\circ\bar{\psi}_i$, which is a $C^{k+1}$-diffeomorpism 
of $M$ onto $M_i$.  
It is easy to  show that there exists an increasing function ${\bf j}:\mathbb N\to\mathbb N$ satisfying 
$\|(\iota_i^{{\bf j}(i)}\circ\psi_i^{{\bf j}(i)})^{\ast}g_{\mathbb E}
-(\bar{\iota}_i^{{\bf j}(i)}\circ\bar{\psi}_i)^{\ast}g_{\mathbb E}\|_{C^k}\to 0\,\,\,\,(i\to\infty)$.  
Easily we can show $\|(\phi_i^{{\bf j}(i)})^{\ast}g_i-g\|_{C^k}\to 0\quad(i\to\infty)$.  
Thus $\{[(M_i,g_i)]\}_{i=1}^{\infty}$ converges to $[(M,g)]$ in the sense of Hamilton.  

Next we shall show the converse.  
Assume that $\{[(M_i,g_i)]\}_{i=1}^{\infty}$ converges to $[(M,g)]$ in the sense of Hamilton.  
Then we can take a sequence $\{\phi_i\}_{i=1}^{\infty}$ of $C^{k+1}$-diffeomorphisms $\phi_i:M\to M_i$'s 
satisfying $\|\phi_i^{\ast}g_i-g\|_{C^k}\to 0$ ($i\to\infty$).  
Take a sequence $\{f^l\}_{l=1}^{\infty}$ in ${\rm Emb}_I((M,g),\mathbb E^m)$ satisfying 

(i) The barycenters of $f^l(M)$'s is equal to the origin $o$ of $\mathbb E^m$;

(ii) $\|f^l\|_{C^k}\to\tau_{\min}^k(M,g)$ ($l\to\infty$).  

\noindent
Then it follows from Lemma 3.2 that, for each $i\in\mathbb N$, take 
a sequence $\{f_i^l\}_{l=1}^{\infty}$ in ${\rm Emb}_I((M_i,g_i),$\newline
$\mathbb E^{m_i})$ having the same barycenter $o$ and 
satisfying the conditions:

(ii') $\|f_i^l\|_{C^k}\to\tau_{\min}^k(M_i,g_i)$ ($l\to\infty$);

(iii) $\displaystyle{\sum_{j=0}^{k+1}\,d_{H,(T^j\mathbb E^{\widehat m},d_{(g_{\mathbb E})_S^j})}
(d^j(\iota_i\circ f_i^l)(S^jM_i),\,d^j(\iota\circ f^l\circ\psi_i)(S^jM))}$

\hspace{5.75truecm} $\displaystyle{\to d_{GH}^k([(M_i,g_i)],\,[(M,g)])\,\,\,\,(l\to\infty)}$\newline

\vspace{-0.15truecm}

for some $\psi_i\in{\rm Diff}^{k+1}(M)$, where $\widehat m$ is a sufficiently large natural number with 
$\widehat m\geq m$ 

and $\widehat m\geq m_i$ ($i\in\mathbb N$), and $\iota$ (resp. $\iota_i$) is suitable totally geodesic 
embedding of $\mathbb E^m$

(resp. $\mathbb E^{m_i}$) into $\mathbb E^{\widehat m}$ which 
maps to the origin of $\mathbb E^m$ (resp. $\mathbb E^{m_i}$) to the origin $o$ 

of $\mathbb E^{\widehat m}$.  

\vspace{0.15truecm}

\noindent
From $\|\phi_i^{\ast}g_i-g\|_{C^k}\to 0$ ($i\to\infty$), we can show 
$$\begin{array}{r}
\displaystyle{\sum_{j=0}^{k+1}\,d_{H,(T^j\mathbb E^{\widehat m_i},d_{(g_{\mathbb E})_S^j})}
(d^j(\iota_i\circ f_i^{{\bf j}(i)})(S^jM_i),\,
d^j(\iota\circ f^{{\bf j}(i)}\circ\psi_i)(S^jM))}\\
\displaystyle{\to 0\quad(i\to\infty)}
\end{array}$$
for some increasing function ${\bf j}:\mathbb N\to\mathbb N$.  Therefore we obtain 
$$d_{GH}^k([(M_i,g_i)],[(M,g)])\to 0\quad(i\to\infty).$$
\qed

\vspace{1truecm}

{\small\textit{Department of Mathematics, Faculty of Science, 
Tokyo University of Science,}}

{\small\textit{1-3 Kagurazaka Shinjuku-ku, Tokyo 162-8601 Japan}}

{\small\textit{E-mail address}: koike@rs.kagu.tus.ac.jp}

\end{document}